\documentclass[12pt,oneside,reqno]{amsart}
\usepackage{mathrsfs}
\usepackage{graphics}
\usepackage{enumerate, amssymb}
\newcommand{\dif}{\mathrm{d}}

\newcommand{\be}{\begin{eqnarray}}
\newcommand{\ee}{\end{eqnarray}}
\newcommand{\ce}{\begin{eqnarray*}}
\newcommand{\de}{\end{eqnarray*}}
\newtheorem{theorem}{Theorem}[section]
\newtheorem{lemma}[theorem]{Lemma}
\newtheorem{remark}[theorem]{Remark}
\newtheorem{definition}[theorem]{Definition}
\newtheorem{proposition}[theorem]{Proposition}
\newtheorem{Example}[theorem]{Example}
\newtheorem{corollary}[theorem]{Corollary}
\def\e{\varepsilon}

\def\a{\alpha}

\def\[{{\Big[}}
\def\]{{\Big]}}
\def\<{{\langle}}
\def\>{{\rangle}}
\def\({{\Big(}}
\def\){{\Big)}}

\def\no{\nonumber}
\def\bt{\begin{theorem}}
\def\et{\end{theorem}}
\def\bl{\begin{lemma}}
\def\el{\end{lemma}}
\def\br{\begin{remark}}
\def\er{\end{remark}}
\def\bx{\begin{Example}}
\def\ex{\end{Example}}
\def\bd{\begin{definition}}
\def\ed{\end{definition}}
\def\bp{\begin{proposition}}
\def\ep{\end{proposition}}
\def\bc{\begin{corollary}}
\def\ec{\end{corollary}}
\def\cA{{\mathcal A}}
\def\cB{{\mathcal B}}
\def\cC{{\mathcal C}}
\def\cD{{\mathcal D}}

\def\cL{{\mathcal L}}
\def\cM{{\mathcal M}}

\def\cP{{\mathcal P}}

\def\mE{{\mathbb E}}

\def\mN{{\mathbb N}}

\def\mP{{\mathbb P}}

\def\mR{{\mathbb R}}

\def\mU{{\mathbb U}}

\def\geq{\geqslant}
\def\leq{\leqslant}

\begin{document}

\allowdisplaybreaks

\title{Convergence of Nonlinear Filterings for Stochastic Dynamical Systems with L\'evy Noises*}

\author{Huijie Qiao}



\thanks{*This work was supported by NSF of China (No. 11001051, 11371352) and the Jiangsu Provincial Key Laboratory of Networked Collective Intelligence under Grant No. BM2017002.}


\date{}

\dedicatory{School of Mathematics \& Jiangsu Provincial Key Laboratory of Networked Collective Intelligence,\\
Southeast University, Nanjing, Jiangsu 211189,  China\\
hjqiaogean@seu.edu.cn}

\begin{abstract}
We consider a nonlinear filtering problem of multiscale non-Gaussian signal processes and observation processes with jumps. Firstly, we prove that the dimension for the signal system can be reduced by a homogenized approach. Secondly, convergence of the corresponding nonlinear filtering to the homogenized filtering is shown by a weak convergence technique. Finally, we give an example to explain our result.
\end{abstract}

\maketitle \rm

\section{Introduction}

Fix $T>0$ and a completed filtered probability space $(\Omega, \mathscr{F}, \{\mathscr{F}_t\}_{t\in[0,T]},\mP)$. 
$V, W$ are $n$-dimensional and $m$-dimensional standard Brownian motion defined on $(\Omega,
\mathscr{F}, \{\mathscr{F}_t\}_{t\in[0,T]},\mP)$, respectively. And $\nu_1, \nu_2$ are two $\sigma$-finite measures defined on a normed space ($\mU,\mathscr{U}, \|\cdot\|_{\mU}$). Let $p_1, p_2$ be two stationary $\mU$-valued Poisson point processes of the class (quasi left-continuous) defined on $(\Omega,\mathscr{F}, \{\mathscr{F}_t\}_{t\in[0,T]},\mP)$ with the characteristic measure $\nu_1, \nu_2$, respectively. Thus, we denote the counting measure of $p_1(t)$ as $N_{p_1}((0,t],\dif u)$ and have $\mE N_{p_1}((0,t],A)=t\nu_1(A)$ for $A\in\mathscr{U}$. Fix $\mU_1, \mU_2\in\mathscr{U}$ with $\nu_1(\mU\setminus\mU_1)<\infty$
and $\nu_2(\mU\setminus\mU_2)<\infty$. Set
\ce
\tilde{N}_{p_1}((0,t], A):=N_{p_1}((0,t], A)-t\nu_1(A), \qquad A\in\mathscr{U}|_{\mU_1},
\de
the compensated martingale measure of $p_1(t)$, and then $\tilde{N}_{p_1}((0,t], A)$ is a $\{\mathscr{F}_t\}_{t\in[0,T]}$-martingale. By the same way, we could define $N_{p_2}((0,t],\dif u), \tilde{N}_{p_2}((0,t],\dif u)$. Besides, we introduce another Poisson random measure $N_{p_2}^{\e}((0,t],\dif u)$ on ($\mU,\mathscr{U}$) such that $\mE N_{p_2}^{\e}((0,t],A)=\frac{1}{\e}t\nu_2(A)$ for $A\in\mathscr{U}$. Moreover, $V_t, W_t, N_{p_1}, N_{p_2}, N_{p_2}^{\e}$ are mutually independent. Consider the following slow-fast system on $\mR^n\times\mR^m$: for $0\leq t\leq T$,
\be\left\{\begin{array}{l}
\dif X^\e_t=b_1(X^\e_t,Z^\e_t)\dif t+\sigma_1(X^\e_t,Z^\e_t)\dif V_t+\int_{\mU_1}f_1(X^\e_{t-}, u)\tilde{N}_{p_1}(\dif t, \dif u), \\
X^\e_0=x_0,\\
\dif Z^\e_t=\frac{1}{\e}b_2(X^\e_t,Z^\e_t)\dif t+\frac{1}{\sqrt{\e}}\sigma_2(X^\e_t,Z^\e_t)\dif W_t+\int_{\mU_2}f_2(X^\e_{t-},Z^\e_{t-},u)\tilde{N}^{\e}_{p_2}(\dif t, \dif u),\\
Z^\e_0=z_0,
\end{array}
\right.
\label{Eq1}
\ee
where these mappings $b_1:\mR^n\times\mR^m\mapsto\mR^n$, $b_2:\mR^n\times\mR^m\mapsto\mR^m$, $\sigma_1:\mR^n\times\mR^m\mapsto\mR^{n\times n}$,
$\sigma_2:\mR^n\times\mR^m\mapsto\mR^{m\times m}$, $f_1:\mR^n\times\mU_1\mapsto\mR^n$ and $f_2:\mR^n\times\mR^m\times\mU_2\mapsto\mR^m$
are all Borel measurable.  

The slow-fast dynamical system (\ref{Eq1}) is usually called multiscale processes, where the rates of change for different variables differ by orders of magnitude. And multiple time scales  models are widely applied in the science and engineering fields (\cite{Dit, kus, kus2, q0, zqd}). For example, fast atmospheric and slow oceanic dynamics describe the climate evolution and state dynamic in electric power systems consists of fast- and slowly-varying elements. 

Next, define an observation process $Y^{\e}$ by
\be
Y_t^{\e}=B_t+\int_0^th(X_s^{\e},Z^\e_s)\dif s+\int_0^t\int_{\mU_3}f_3(s,u)\tilde{N}_{\lambda}(\dif s, \dif u)+\int_0^t\int_{\mU\setminus\mU_3}g_3(s,u)N_{\lambda}(\dif s, \dif u),
\label{ypde}
\ee
where $B$ is a $d$-dimensional standard Brownian motion and $N_{\lambda}(\dif t,\dif u)$ is a random measure with a predictable compensator $\lambda(t,X_t^\e,u)\dif t\nu_3(\dif u)$. Here the function $\lambda: [0,T]\times\mR^n\times\mU\rightarrow(0,1)$ is Borel measurable and $\nu_3$ is a $\sigma$-finite measure defined on $\mU$ with $\nu_3(\mU\setminus\mU_3)<\infty$ and $\int_{\mU_3}\|u\|_{\mU}^2\,\nu_3(\dif u)<\infty$ for a fixed $\mU_3\in\mathscr{U}$. Concretely speaking, set 
$$
\tilde{N}_\lambda((0,t], A):=N_\lambda((0,t],A)-\int_0^t\int_A\lambda(s,X_s^\e,u)\dif s\nu_3(\dif u), \quad t\in[0,T], A\in\mathscr{U}|_{\mU_3},
$$ 
and then $\tilde{N}_\lambda((0,t],\dif u)$ is the compensated martingale measure
of $N_{\lambda}((0,t],\dif u)$. Moreover, $V_t, W_t, B_t, N_{p_1}, N_{p_2}, N_{p_2}^{\e}, N_{\lambda}$ are mutually independent. These functions $h: \mR^n\times\mR^m\mapsto\mR^d$, $f_3:[0,T]\times\mU_3\mapsto\mR^d$ and $g_3:[0,T]\times(\mU\setminus\mU_3)\mapsto\mR^d$
are all Borel measurable. For a Borel measurable function $F$, the nonlinear filtering  problem for the slow component $X_t^{\e}$ with respect to $\{Y^{\e}_s, 0\leq s\leq t\}$ leads to evaluating the `filter' $\mE[F(X^{\e}_t)|\mathscr{F}_t^{Y^{\e}}]$, where $\mathscr{F}_t^{Y^{\e}}$ is the $\sigma$-algebra generated by $\{Y^{\e}_s, 0\leq s\leq t\}$ and $\mE|F(X^{\e}_t)|<\infty$ for $t\in[0,T]$.

When $f_1=f_2=f_3=g_3=0$, this problem has been studied alternatively. Let us recall some works. 
In \cite{psn}, Park-Sowers-Namachchivaya considered the filtering problem with a two-dimensional plant and a one-dimensional observation process. There they used the time change and decomposition methods. And then for the high dimensional case, Park-Namachchivaya-Yeong \cite{pny} presented a numerical algorithm method. Later, Imkeller-Namachchivaya-Perkowski-Yeong \cite{ImkellerSri} showed that for the high dimensional slow-fast dynamical system (\ref{Eq1}), the filter  $\mE[F(X^{\e}_t)|\mathscr{F}_t^{Y^{\e}}]$ converges to the homogenized filter (See Section \ref{filpro}) by double backward stochastic differential equations and asymptotic techniques. Recently, Papanicolaou-Spiliopoulos \cite{ps1} also studied this convergence problem by independent version technique  and then applied it to statistical inference. When $f_1\neq0, f_2\neq0$, $f_3=g_3=0$, Kushner \cite{kus} studied this problem by a weak convergence method.

Recently, some filtering problems, whose observation parts are driven by It\^o-L\'evy processes,  appear in geophysics and finance (\cite{Dit, mbfp}). For example, $Y_t^{\e}$ stands for an observed asset price process in \cite{mbfp}. Moreover, some data information is well described only by L\'evy processes or It\^o-L\'evy processes. That prompts us to notice the case with $f_1\neq0, f_2\neq0, f_3\neq0$, or $g_3\neq0$. First of all, the dimension for the slow-fast system is proved to be reduced by a homogenized approach. And then, we show that the nonlinear filtering of the origin slow-fast system also converges to the homogenized filtering.

It is worthwhile to mention our methods. For the dimension reduction of the slow-fast systems, at present there are two methods: homogenization (\cite{ImkellerSri}) and invariant manifolds (\cite{q0, zqd}). The former heavily depends on martingale problems for origin systems and the latter needs good dynamical structures. Therefore, based on the form of the system (\ref{Eq1}), we choose the first method to reduce it. Next, for the filtering problem of the slow component $X_t^{\e}$ with respect to $\{Y^{\e}_s, 0\leq s\leq t\}$, since the time change is only useful for a one-dimensional process, the technique can {\it not} be applied  to the present case. Moreover, the theory for double backward stochastic differential equations with jumps is short. Therefore, we compute the difference between $\mE[F(X^{\e}_t)|\mathscr{F}_t^{Y^{\e}}]$ and the homogenized filter and then convert it to the difference between two unnormalized filterings. With the help of the weak convergence method in \cite{kus}, we know that the difference between two unnormalized filterings converges to zero. Thus, we prove that $\mE[F(X^{\e}_t)|\mathscr{F}_t^{Y^{\e}}]$ converges weakly to the homogenized filter. Although our method is the same to one in \cite{kus}, some new techniques are used due to jumps, such as ones in constructing two exponential martingales and deducing two equations which two unnormalized filterings satisfy.

The paper is arranged as follows. In the next section, we introduce some notation, terminology and concepts used in the sequel. The dimensional reduction for these slow-fast systems is placed in Section \ref{conpro}. In Section \ref{filpro}, nonlinear filtering problems are introduced. And convergence of the nonlinear filtering for the origin slow-fast system to the homogenized filtering is proved in Section \ref{confil}. In Section \ref{exam} we give an example to explain our result. We summarize the paper and mention some future questions in Section \ref{con}.

The following convention will be used throughout the paper: $C$ with
or without indices will denote different positive constants
(depending on the indices) whose values may change from one place to
another.

\section{Preliminary}\label{pre}

In the section, we introduce some notation, terminology, concepts and state a theorem used in the sequel.
Firstly, introduce the following notation and terminology:

$(i)$ For a separable metric space $E$, let $\mathscr{B}(E)$ denote the Borel $\sigma$-algebra on $E$. Let $\cB(E), \cC(E)$ denote the set of all real-valued uniformly bounded $\mathscr{B}(E)$-measurable and continuous functions on $E$, respectively. Let $\bar{\cC}(E) := \cB(E)\cap \cC(E)$. Let $\cC_c(E)$ be the
set of all members of $\bar{\cC}(E)$ which have compact supports. When $E$ is locally compact, let $\hat{\cC}(E)$ be the collection of all
members of $\bar{\cC}(E)$ which vanish at infinity.

$(ii)$ For a positive integer $r$, let $\cC^r(\mR^q)$ denote the collection of all members of $\cC(\mR^q)$ with
continuous derivatives of each order, up to and including $r$. Let $\cC_c^\infty(\mR^q)$ denote the collection of all members of
$\cC(\mR^q)$ with continuous derivatives of all orders and compact supports. For a metric space $E$  and some
positive integer $r$, $\cC^{r, 0}(\mR^q\times E)$ denotes the collection of all mappings $f\in \cC(\mR^q\times E)$
whose partial derivatives of every order up to and including $r$, with respect to its first $q$
real-valued arguments, exist and are members of $C(\mR^q\times E)$. Let $C_c^{r, 0}(\mR^q\times E):=
C^{r, 0}(\mR^q\times E)\cap C_c(\mR^q\times E)$.

$(iii)$ For a complete separable metric space $E$, let $\cP(E)$ denote the collection of
all probability measures on the measurable space $(E, \mathscr{B}(E))$ with the usual topology of
weak convergence. Let $L(X)$ be the distribution of $E$-valued random variable $X$. Set $\mu f := \int_E f \dif\mu$ for $f\in \cB(E)$.

Secondly, we introduce some concepts. Suppose that $E$ is a separable metric space.

\bd\label{martprob}
Let $\cA\subset B(E)\times B(E)$ be a relation with the domain $\cD(\cA)$, and let $\mu\in\cP(E)$.
Then a progressively measurable solution of the martingale problem for $\cA$ (for $(\cA,\mu)$) is some pair
$\{(\tilde{\Omega}, \tilde{\mathscr{F}}, \{\tilde{\mathscr{F}}_t\}_{t\in[0,T]}, \tilde{\mP}), (\tilde{X}_t)\}$, in which $(\tilde{\Omega},
\tilde{\mathscr{F}}, \{\tilde{\mathscr{F}}_t\}_{t\in[0,T]}, \tilde{\mP})$ is a completed filtered probability space and $\{\tilde{X}_t\}$
is an $E$-valued $\{\tilde{\mathscr{F}}_t\}_{t\in[0,T]}$-progressively measurable process such that $f(\tilde{X}_t)-\int_0^t\cA f(\tilde{X}_s)\dif s$
is an $\{\tilde{\mathscr{F}}_t\}_{t\in[0,T]}$-martingale for each $f\in\cD(\cA)$ (and $L(\tilde{X}_0)=\mu$). The martingale problem
for $(\cA, \mu)$ has the property of existence when there exists some progressively measurable solution
of the martingale problem for $(\cA, \mu)$, and has the property of uniqueness when, given any two progressively
measurable solutions $\{(\tilde{\Omega}, \tilde{\mathscr{F}}, \{\tilde{\mathscr{F}}_t\}_{t\in[0,T]}, \tilde{\mP}), (\tilde{X}_t)\}$ and
$\{(\check{\Omega}, \check{\mathscr{F}}, \{\check{\mathscr{F}}_t\}_{t\in[0,T]}, \check{\mP}), (\check{X}_t)\}$ of the martingale problem
for $(\cA, \mu)$, the $E$-valued processes $\tilde{X}$ and $\check{X}$ necessarily have identical
finite-dimensional distributions. The martingale problem for $(\cA,\mu)$ is called well-posed when it has the
properties of both existence and uniqueness. Finally, the martingale problem for $\cA$ is well-posed
when the martingale problem for $(\cA, \mu)$ is well-posed for each $\mu\in\cP(E)$.
\ed

Thirdly, the following theorem comes from Theorem 2.7 in \cite{kur}.

\bt\label{tigth}
Let $\{X_n\}$ be a sequence of stochastic processes with sample paths in $D([0,\infty),E)$, where $E$ is a complete and separable metric space with the metric $d$. Let $\mathscr{F}_t^n=\sigma(X_n(s): s\leq t)$ and $M_T^n$ be the collection of $(\mathscr{F}_t^n)$-stopping times $\tau$ such that $\tau\leq T$ a.s.. Suppose 

(i) for every $\eta>0$ and rational $t\geq 0$ there exists a compact set $\Gamma_{t, \eta}\subset E$ such that 
$$
\inf_n\mP\{X_n(t)\in\Gamma_{t, \eta}\}>1-\eta;
$$

(ii) for every $T>0$ and $\beta>0$
$$
\lim_{\delta\rightarrow0}\limsup_{n\rightarrow\infty}\sup_{\tau\in M_T^n}E[d^\beta(X_n(\tau),X_n(\tau+\delta))]=0.
$$

Then $\{X_n\}$ is relatively weakly compact.
\et

\section{Convergence of the system (\ref{Eq1})}\label{conpro}

In the section, we study convergence for the system (\ref{Eq1}) when $\e\rightarrow0$.

\subsection{The fast-slow equation (\ref{Eq1})}\label{fasl}

In the subsection, we consider the system (\ref{Eq1}), i.e. for $0\leq t\leq T$,
\ce\left\{\begin{array}{l}
\dif X^\e_t=b_1(X^\e_t,Z^\e_t)\dif t+\sigma_1(X^\e_t,Z^\e_t)\dif V_t+\int_{\mU_1}f_1(X^\e_{t-}, u)\tilde{N}_{p_1}(\dif t, \dif u), \\
X^\e_0=x_0,\\
\dif Z^\e_t=\frac{1}{\e}b_2(X^\e_t,Z^\e_t)\dif t+\frac{1}{\sqrt{\e}}\sigma_2(X^\e_t,Z^\e_t)\dif W_t+\int_{\mU_2}f_2(X^\e_{t-},Z^\e_{t-},u)\tilde{N}^{\e}_{p_2}(\dif t, \dif u),\\
Z^\e_0=z_0.
\end{array}
\right.
\de
We give out our assumptions and state some related results.

\begin{enumerate}[\bf{Assumption 1.}]
\item
\end{enumerate}
\begin{enumerate}[($\mathbf{H}^1_{b_1, \sigma_1, f_1}$)]
\item For $x_1, x_2\in\mR^n$, $z_1, z_2\in\mR^m$, there exist $L_{b_1}, L_{\sigma_1}, L_{f_1}>0$ such that
\ce
&|b_1(x_1, z_1)-b_1(x_2, z_2)|\leq L_{b_1}(|x_1-x_2|+|z_1-z_2|),\\
&\|\sigma_1(x_1, z_1)-\sigma_1(x_2, z_2)\|\leq L_{\sigma_1}(|x_1-x_2|+|z_1-z_2|),\\
&\int_{\mU_1}|f_1(x_1,u)-f_1(x_2,u)|^2\,\nu_1(\dif u)\leq L_{f_1}|x_1-x_2|^2,
\de
where $|\cdot|$ and $\|\cdot\|$ denote the length of a vector and the Hilbert-Schmidt norm of
a matrix, respectively.
\end{enumerate}

\begin{enumerate}[($\mathbf{H}^2_{b_1, \sigma_1, f_1}$)]
\item For $x\in\mR^n$, $z\in\mR^m$, there exists a $L_{b_1, \sigma_1, f_1}>0$ such that
$$
|b_1(x,z)|^2+\|\sigma_1(x,z)\|^2+\int_{\mU_1}|f_1(x,u)|^2\nu_1(\dif u)\leq L_{b_1, \sigma_1, f_1}.
$$
\end{enumerate}

\begin{enumerate}[($\mathbf{H}^1_{b_2}$)]
\item (i) $b_2$ is bi-continuous in $(x, z)$,\\
(ii) There exist $L_{b_2}\geq0, \bar{L}_{b_2}>0$ such that
\ce
&&|b_2(x_1, z)-b_2(x_2, z)|\leq L_{b_2}|x_1-x_2|, \qquad\qquad\qquad x_1, x_2\in\mR^n, z\in\mR^m,\\
&&\<z_1-z_2, b_2(x, z_1)-b_2(x, z_2)\>\leq -\bar{L}_{b_2}|z_1-z_2|^2, \qquad x\in\mR^n, z_1, z_2\in\mR^m,
\de
(iii) For $x\in\mR^n$, $z\in\mR^m$, there exists a constant $\bar{\bar{L}}_{b_2}>0$ such that
$$
|b_2(x,z)|\leq \bar{\bar{L}}_{b_2}(1+|x|+|z|).
$$
\end{enumerate}

\begin{enumerate}[($\mathbf{H}^1_{\sigma_2}$)]
\item For $x_1, x_2\in\mR^n$, $z_1, z_2\in\mR^m$, there exists a constant $L_{\sigma_2}>0$ such that
\ce
\|\sigma_2(x_1, z_1)-\sigma_2(x_2, z_2)\|^2\leq L_{\sigma_2}(|x_1-x_2|^2+|z_1-z_2|^2).
\de
\end{enumerate}

\begin{enumerate}[($\mathbf{H}^1_{f_2}$)]
\item There exists a positive function $L(u)$ satisfying 
\ce 
\sup_{u\in\mU_2}L(u)\leq\gamma<1~\mbox{ and } \int_{\mU_2}L(u)^2\,\nu_2(\dif u)<+\infty, 
\de 
such that for any $x_1, x_2\in\mR^n$, $z_1, z_2\in\mR^m$ and $u\in\mU_2$ 
\ce
|f_2(x_1,z_1,u)-f_2(x_2,z_2,u)|\leq L(u)(|x_1-x_2|+|z_1-z_2|),
\de 
and 
\ce 
|f_2(0,0,u)|\leq L(u). 
\de
\end{enumerate}

\medspace

Under {\bf Assumption 1.}, by Theorem 1.2 in \cite{q2}, the system (\ref{Eq1}) has a unique strong solution denoted
by $(X^\e_t,Z^\e_t)$. Moreover, the infinitesimal generator of the system (\ref{Eq1}) is given by
\ce
(\cL^{\e}H)(x,z)=(\cL^{X^\e}H)(x,z)+(\cL^{Z^\e}H)(x,z), \qquad H\in\cD(\cL^{\e}),
\de
where
\ce
(\cL^{X^\e}H)(x,z)&:=&\frac{\partial H(x,z)}{\partial x_i}b^i_1(x,z)+\frac{1}{2}\frac{\partial^2H(x,z)}{\partial x_i\partial x_j}
(\sigma_1\sigma_1^T)^{ij}(x,z)\\
&&+\int_{\mU_1}\Big[H\big(x+f_1(x,u),z\big)-H(x,z)
-\frac{\partial H(x,z)}{\partial x_i}f^i_1(x,u)\Big]\nu_1(\dif u),
\de
and
\ce
(\cL^{Z^\e}H)(x,z)&:=&\frac{1}{\e}\frac{\partial H(x,z)}{\partial z_i}b^i_2(x,z)+\frac{1}{2\e}\frac{\partial^2H(x,z)}{\partial z_i\partial z_j}
(\sigma_2\sigma_2^T)^{ij}(x,z)\\
&&+\frac{1}{\e}\int_{\mU_2}\Big[H\big(x,z+f_2(x,z,u)\big)-H(x,z)
-\frac{\partial H(x,z)}{\partial z_i}f^i_2(x,z,u)\Big]\nu_2(\dif u).
\de
Here and hereafter, we use the convention that repeated indices imply summation.

\subsection{The fast equation}\label{fas}

In the subsection, we mainly study the second part of the system (\ref{Eq1}). 

First, take any $x\in\mR^n$ and fix it. And consider the following SDE in $\mR^m$:
\ce\left\{\begin{array}{l}
\dif Z^x_t=b_2(x,Z^x_t)\dif t+\sigma_2(x,Z^x_t)\dif W_t+\int_{\mU_2}f_2(x,Z^x_t,u)\tilde{N}_{p_2}(\dif t, \dif u),\\
Z^x_0=z_0, \qquad t\geq0.
\end{array}
\right.
\de
Under the assumption ($\mathbf{H}^1_{b_2}$) ($\mathbf{H}^1_{\sigma_2}$) ($\mathbf{H}^1_{f_2}$), the above equation has a unique solution $Z^x_t$.
In addition, it is a Markov process and its transition probability is denoted by $p(x; z_0,t,A)$
for $t\geq0$ and $A\in\mathscr{B}(\mR^m)$. Set $(T^x_t\varphi)(z_0):=\int_{\mR^m}\varphi(z')p(x; z_0,t,\dif z')$ for any
$\varphi\in \cC(\mR^m)$, and then $\{T^x_t, t\geq0\}$ is its transition semigroup and $\e\cL^{Z^\e}$ is its infinitesimal generator. For $Z^x_t$, we assume:

\begin{enumerate}[\bf{Assumption 2.}]
\item
\end{enumerate}
\begin{enumerate}[($\mathbf{H}^2_{\sigma_2}$)]
\item There exists a strictly positive function $\a_1(x)$ such that
\ce
\<\sigma_2(x,z)h,h\>\geq\sqrt{\a_1(x)}|h|^2, \qquad z,h\in\mR^m,
\de
and 
\ce
\|\sigma_{\a_1}(x,z_1)-\sigma_{\a_1}(x,z_2)\|^2\leq L_{\a_1}|z_1-z_2|^2, \qquad z_1, z_2\in\mR^m,
\de 
where $\sigma_{\a_1}(x,z)$ is the unique symmetric nonnegative definite matrix-valued function
such that $\sigma_{\a_1}(x,z)\sigma_{\a_1}(x,z)
=\sigma_2(x,z)\sigma^T_2(x,z)-\a_1(x)\emph{I}$ for the unit matrix $\emph{I}$.
\end{enumerate}
\begin{enumerate}[($\mathbf{H}^1_{b_2,\sigma_2,f_2}$)] 
\item There exist a $r\geq 2$ and two functions $\a_2(x)>0$, $\a_3(x)\geq0$ such that for all $z\in\mR^m$ 
\ce
2\<z,b_2(x,z)\>+\|\sigma_2(x,z)\|^2+\int_{\mU_2}\big|f_2(x,z,u)\big|^2\nu_2(\dif
u)\leq-\a_2(x)|z|^r+\a_3(x). 
\de
\end{enumerate}
\begin{enumerate}[($\mathbf{H}^2_{b_2,\sigma_2,f_2}$)] 
\item $$
M:=2\bar{L}_{b_2}-L_{b_2}-L_{\sigma_2}-2\int_{\mU_2}L^2(u)\nu_2(\dif u)>0
$$
\end{enumerate}

Under the assumptions ($\mathbf{H}^1_{b_2}$) ($\mathbf{H}^1_{\sigma_2}$) ($\mathbf{H}^1_{f_2}$) ($\mathbf{H}^2_{\sigma_2}$)($\mathbf{H}^1_{b_2,\sigma_2,f_2}$), by Theorem 1.3 in  \cite{q1} there exists a unique invariant probability measure $\bar{p}(x,\cdot)$  for $Z^x_t$. And set
\ce
\bar{b}_1(x):=\int_{\mR^m}b_1(x,z)\bar{p}(x,\dif z), 
\de
and then we have the following result.

\bl\label{bblico}
Under the assumptions ($\mathbf{H}^1_{b_1, \sigma_1, f_1}$) ($\mathbf{H}^1_{b_2}$) ($\mathbf{H}^1_{\sigma_2}$) ($\mathbf{H}^1_{f_2}$) ($\mathbf{H}^2_{b_2,\sigma_2,f_2}$), there exists a constant $L_{\bar{b}_1}\geq0$ such that for $x_1, x_2\in\mR^n$
\ce
|\bar{b}_1(x_1)-\bar{b}_1(x_2)|\leq L_{\bar{b}_1}|x_1-x_2|.
\de
\el
\begin{proof}
First, by the definition of $\bar{b}_1$ and ($\mathbf{H}^1_{b_1, \sigma_1, f_1}$) it holds that
\be
|\bar{b}_1(x_1)-\bar{b}_1(x_2)|&=&\left|\lim_{S\rightarrow\infty}\frac{1}{S}\int_0^S\mE b_1(x_1,Z_t^{x_1})\dif t-\lim_{S\rightarrow\infty}\frac{1}{S}\int_0^S\mE b_1(x_2,Z_t^{x_2})\dif t\right|\no\\
&=&\left|\lim_{S\rightarrow\infty}\frac{1}{S}\int_0^S\mE(b_1(x_1,Z_t^{x_1})-b_1(x_2,Z_t^{x_2}))\dif t\right|\no\\
&\leq&\lim_{S\rightarrow\infty}\frac{1}{S}\int_0^S\mE\left|b_1(x_1,Z_t^{x_1})-b_1(x_2,Z_t^{x_2})\right|\dif t\no\\
&\leq&L_{b_1}|x_1-x_2|+L_{b_1}\lim_{S\rightarrow\infty}\frac{1}{S}\int_0^S\mE\left|Z_t^{x_1}-Z_t^{x_2}\right|\dif t\no\\
&\leq&L_{b_1}|x_1-x_2|+L_{b_1}\lim_{S\rightarrow\infty}\(\frac{1}{S}\int_0^S\mE\left|Z_t^{x_1}-Z_t^{x_2}\right|^2\dif t\)^{1/2},
\label{bb1}
\ee
where the last inequality is based on the H\"older inequality. And then we estimate $\frac{1}{S}\int_0^S\mE\left|Z_t^{x_1}-Z_t^{x_2}\right|^2\dif t$. For $Z_t^{x_1}, Z_t^{x_2}$, it follows from the It\^o formula that
\ce
\left|Z_t^{x_1}-Z_t^{x_2}\right|^2&=&2\int_0^t\<Z_s^{x_1}-Z_s^{x_2}, b_2(x_1,Z_s^{x_1})-b_2(x_2,Z_s^{x_2})\>\dif s\\
&&+2\int_0^t\<Z_s^{x_1}-Z_s^{x_2}, (\sigma_2(x_1,Z_s^{x_1})-\sigma_2(x_2,Z_s^{x_2}))\dif W_s\>\\
&&+\int_0^t\int_{\mU_2}[|Z_s^{x_1}-Z_s^{x_2}+f_2(x_1,Z^{x_1}_s,u)-f_2(x_2,Z^{x_2}_s,u)|^2\\
&&\qquad -|Z_s^{x_1}-Z_s^{x_2}|^2]\tilde{N}_{p_2}(\dif s, \dif u)\\
&&+\int_0^t\|\sigma_2(x_1,Z_s^{x_1})-\sigma_2(x_2,Z_s^{x_2})\|^2\dif s\\
&&+\int_0^t\int_{\mU_2}|f_2(x_1,Z^{x_1}_s,u)-f_2(x_2,Z^{x_2}_s,u)|^2\nu_2(\dif u)\dif s.
\de
Taking the expectation on two sides, one can obtain 
\be
\mE\left|Z_t^{x_1}-Z_t^{x_2}\right|^2&=&2\mE\int_0^t\<Z_s^{x_1}-Z_s^{x_2}, b_2(x_1,Z_s^{x_1})-b_2(x_2,Z_s^{x_2})\>\dif s\no\\
&&+\int_0^t\mE\|\sigma_2(x_1,Z_s^{x_1})-\sigma_2(x_2,Z_s^{x_2})\|^2\dif s\no\\
&&+\int_0^t\int_{\mU_2}\mE|f_2(x_1,Z^{x_1}_s,u)-f_2(x_2,Z^{x_2}_s,u)|^2\nu_2(\dif u)\dif s\no\\
&=:&I_1+I_2+I_3.
\label{i123}
\ee
For $I_1$, by ($\mathbf{H}^1_{b_2}$) it holds that
\be
I_1&=&2\mE\int_0^t\<Z_s^{x_1}-Z_s^{x_2}, b_2(x_1,Z_s^{x_1})-b_2(x_1,Z_s^{x_2})\>\dif s\no\\
&&+2\mE\int_0^t\<Z_s^{x_1}-Z_s^{x_2}, b_2(x_1,Z_s^{x_2})-b_2(x_2,Z_s^{x_2})\>\dif s\no\\
&\leq&-2\bar{L}_{b_2}\int_0^t\mE|Z_s^{x_1}-Z_s^{x_2}|^2\dif s+2L_{b_2}\mE\int_0^t|Z_s^{x_1}-Z_s^{x_2}||x_1-x_2|\dif s\no\\
&\leq&(-2\bar{L}_{b_2}+L_{b_2})\int_0^t\mE|Z_s^{x_1}-Z_s^{x_2}|^2\dif s+L_{b_2}|x_1-x_2|^2 t.
\label{i1}
\ee
For $I_2$, by ($\mathbf{H}^1_{\sigma_2}$) we know that
\be
I_2&\leq&L_{\sigma_2}\int_0^t\mE(|x_1-x_2|^2+|Z_s^{x_1}-Z_s^{x_2}|^2)\dif s\no\\
&=&L_{\sigma_2}\int_0^t\mE|Z_s^{x_1}-Z_s^{x_2}|^2\dif s+L_{\sigma_2}|x_1-x_2|^2 t.
\label{i2}
\ee
By the same deduction to that for $I_2$, it follows from ($\mathbf{H}^1_{f_2}$) that 
\be
I_3&\leq&2\int_{\mU_2}L^2(u)\nu_2(\dif u)\int_0^t\mE|Z_s^{x_1}-Z_s^{x_2}|^2\dif s+2\int_{\mU_2}L^2(u)\nu_2(\dif u)|x_1-x_2|^2 t.
\label{i3}
\ee
Combining (\ref{i1}) (\ref{i2}) (\ref{i3}) with (\ref{i123}), we get that
\ce
\mE\left|Z_t^{x_1}-Z_t^{x_2}\right|^2\leq -M\int_0^t\mE|Z_s^{x_1}-Z_s^{x_2}|^2\dif s+\(L_{b_2}+L_{\sigma_2}+2\int_{\mU_2}L^2(u)\nu_2(\dif u)\)|x_1-x_2|^2 t.
\de
Thus, we have furthermore 
\be
\frac{1}{S}\int_0^S\mE|Z_t^{x_1}-Z_t^{x_2}|^2\dif t\leq \frac{L_{b_2}+L_{\sigma_2}+2\int_{\mU_2}L^2(u)\nu_2(\dif u)}{M}|x_1-x_2|^2.
\label{zes}
\ee
That is, it holds that
\ce
|\bar{b}_1(x_1)-\bar{b}_1(x_2)|\leq L_{\bar{b}_1}|x_1-x_2|,
\de
where $L_{\bar{b}_1}:=L_{b_1}\left(1+\sqrt{\frac{L_{b_2}+L_{\sigma_2}+2\int_{\mU_2}L^2(u)\nu_2(\dif u)}{M}}\right)$. The proof is completed.
\end{proof}

Next, set
\ce
\Sigma(x):=\int_{\mR^m}(\sigma_1\sigma^T_1)(x,z)\bar{p}(x,\dif z),
\de
and then we know that $\Sigma(x)$ is a nonnegative definite symmetric matrix. 

\begin{enumerate}[\bf{Assumption 3.}]
\item
\end{enumerate}
\begin{enumerate}[]
\item  There exists a constant $l>0$ such that
$$
\<\sigma_1(x,z)h,h\>\geq l|h|^2, \quad x,h\in\mR^n, z\in\mR^m.
$$
\end{enumerate}

Thus, based on {\bf Assumption 3.}, it holds that $\Sigma(x)$ is positive definite and $\|\Sigma(x)\|\geq l^2$. Moreover, there exists a unique positive definite symmetric matrix $\bar{\sigma}_1(x)$ such that $\Sigma(x)=(\bar{\sigma}_1\bar{\sigma}_1)(x)$ and $\|\bar{\sigma}_1(x)\|\geq l$. Besides, $\bar{\sigma}_1(x)$ has the following property.

\bl\label{siglip}
Suppose that ($\mathbf{H}^1_{b_1, \sigma_1, f_1}$) ($\mathbf{H}^1_{b_2}$) ($\mathbf{H}^1_{\sigma_2}$) ($\mathbf{H}^1_{f_2}$) ($\mathbf{H}^2_{b_2,\sigma_2,f_2}$) and {\bf Assumption 3.} hold. Then there exists a constant $L_{\bar{\sigma}_1}$ such that for $x_1, x_2\in\mR^n$
\ce
\|\bar{\sigma}_1(x_1)-\bar{\sigma}_1(x_2)\|\leq L_{\bar{\sigma}_1}|x_1-x_2|.
\de
\el
\begin{proof}
First of all, we observe $\Sigma(x)$. By the definition of $\Sigma(x)$, it holds that for $x_1, x_2\in\mR^n$,
\ce
\|\Sigma(x_1)-\Sigma(x_2)\|&=&\|\int_{\mR^m}(\sigma_1\sigma^T_1)(x_1,z)\bar{p}(x_1,\dif z)-\int_{\mR^m}(\sigma_1\sigma^T_1)(x_2,z)\bar{p}(x_2,\dif z)\|\\
&=&\|\lim_{S\rightarrow\infty}\frac{1}{S}\int_0^S\mE(\sigma_1\sigma^T_1)(x_1,Z_t^{x_1})\dif t-\lim_{S\rightarrow\infty}\frac{1}{S}\int_0^S\mE(\sigma_1\sigma^T_1)(x_2,Z_t^{x_2})\dif t\|\\
&\leq&\lim_{S\rightarrow\infty}\frac{1}{S}\int_0^S\mE\|(\sigma_1\sigma^T_1)(x_1,Z_t^{x_1})-(\sigma_1\sigma^T_1)(x_2,Z_t^{x_2})\|\dif t\\
&\leq&\lim_{S\rightarrow\infty}\frac{1}{S}\int_0^S\mE\|\sigma_1(x_1,Z_t^{x_1})\sigma^T_1(x_1,Z_t^{x_1})-\sigma_1(x_2,Z_t^{x_2})\sigma^T_1(x_1,Z_t^{x_1})\|\dif t\\
&&+\lim_{S\rightarrow\infty}\frac{1}{S}\int_0^S\mE\|\sigma_1(x_2,Z_t^{x_2})\sigma^T_1(x_1,Z_t^{x_1})-\sigma_1(x_2,Z_t^{x_2})\sigma^T_1(x_2,Z_t^{x_2})\|\dif t\\
&\leq&2L_{\sigma_1}L^{1/2}_{b_1,\sigma_1,f_1}|x_1-x_2|+2L_{\sigma_1}L^{1/2}_{b_1,\sigma_1,f_1}\lim_{S\rightarrow\infty}\frac{1}{S}\int_0^S\mE|Z_t^{x_1}-Z_t^{x_2}|\dif t\\
&\leq&C|x_1-x_2|,
\de
where ($\mathbf{H}^1_{b_1, \sigma_1, f_1}$)($\mathbf{H}^2_{b_1, \sigma_1, f_1}$) and (\ref{zes}) are used. Note that 
$$
(\bar{\sigma}_1(x_1)-\bar{\sigma}_1(x_2))(\bar{\sigma}_1(x_1)+\bar{\sigma}_1(x_2))=\Sigma(x_1)-\Sigma(x_2).
$$
So, we obtain that
\ce
\|\bar{\sigma}_1(x_1)-\bar{\sigma}_1(x_2)\|&=&\|(\Sigma(x_1)-\Sigma(x_2))(\bar{\sigma}_1(x_1)+\bar{\sigma}_1(x_2))^{-1}\|\\
&\leq&\|\Sigma(x_1)-\Sigma(x_2)\|\|(\bar{\sigma}_1(x_1)+\bar{\sigma}_1(x_2))^{-1}\|\\
&\leq&C|x_1-x_2|.
\de
The proof is over.
\end{proof}

\subsection{The average equation}\label{ave}

Next, we construct a SDE on the probability space $(\Omega, \mathscr{F}, \{\mathscr{F}_t\}_{t\in[0,T]}, \mP)$ as follows:
\be\left\{\begin{array}{l}
\dif X^0_t=\bar{b}_1(X^0_t)\dif t+\bar{\sigma}_1(X^0_t)\dif V_t+\int_{\mU_1}f_1(X^0_{t-}, u)\tilde{N}_{p_1}(\dif t, \dif u),\\
X^0_0=x_0, \qquad\qquad 0\leq t\leq T.
\label{appequ}
\end{array}
\right.
\ee
By Lemma \ref{bblico} and \ref{siglip}, it holds that Eq.(\ref{appequ}) has a unique strong solution denoted as $X^0_t$. Moreover, it is a Markov process. Define an operator $\bar{\cL}$ as follows:
\ce
\cD(\bar{\cL})&:=&\cC_c^\infty(\mR^n),\\
(\bar{\cL}g)(x)&:=&\frac{\partial g(x)}{\partial x_i}\bar{b}^i_1(x)+\frac{1}{2}\frac{\partial^2g(x)}{\partial x_i\partial x_j}
(\bar{\sigma}_1\bar{\sigma}_1)^{ij}(x)\\
&&+\int_{\mU_1}\Big[g\big(x+f_1(x,u)\big)-g(x)
-\frac{\partial g(x)}{\partial x_i}f^i_1(x,u)\Big]\nu_1(\dif u), \quad g\in\cD(\bar{\cL}).
\de
Thus, $\bar{\cL}$ is the infinitesimal generator of $X^0_{\cdot}$. And it follows from \cite{st} that the martingale problem for ($\bar{\cL}$, $\delta_{x_0}$) is well-posed. Next, we state and prove the main result in the section.

\bt\label{conprot}
Under all the above hypotheses $\{X_t^\e, t\in[0,T]\}$ converges weakly to $\{X^0_t, t\in[0,T]\}$ in $D([0,T],\mR^n)$.
\et
\begin{proof}
{\bf Step 1.} We prove that $\{X_t^\e, t\in[0,T]\}$ is relatively weakly compact in $D([0,T],\mR^n)$.

First of all, consider the martingale problem associated with $\cL^{\e}$. For $H\in\cD(\cL^{\e})$,
\be
M_H(t):=H(X^\e_t,Z^\e_t)-H(x_0,z_0)-\int_0^t(\cL^{\e}H)(X^\e_s,Z^\e_s)\dif s
\label{xzmartprob}
\ee
is a square integrable martingale and
\ce
\<M_H(\cdot), M_H(\cdot)\>_t&=&\int_0^t\frac{\partial H(X^\e_s,Z^\e_s)}{\partial x_i}\frac{\partial H(X^\e_s,Z^\e_s)}{\partial x_j}
(\sigma_1\sigma_1^T)^{ij}(X^\e_s,Z^\e_s)\dif s\\
&&+\int_0^t\int_{\mU_1}[H(X^\e_s+f_1(X^\e_s,u),Z^\e_s)-H(X^\e_s,Z^\e_s)]^2\nu_1(\dif u)\dif s\\
&&+\frac{1}{\e}\int_0^t\frac{\partial H(X^\e_s,Z^\e_s)}{\partial z_i}\frac{\partial H(X^\e_s,Z^\e_s)}{\partial z_j}
(\sigma_2\sigma_2^T)^{ij}(X^\e_s,Z^\e_s)\dif s\\
&&+\frac{1}{\e}\int_0^t\int_{\mU_2}[H(X^\e_s,Z^\e_s+f_2(X^\e_s, Z^\e_s, u))-H(X^\e_s,Z^\e_s)]^2\nu_2(\dif u)\dif s.
\de

Taking $H(x,z)=x^i, i=1, 2, \cdots, n$ in (\ref{xzmartprob}), we obtain that
\ce
(X^\e_t)^i-x_0^i=\int_0^tb^i_1(X^\e_r,Z^\e_r)\dif r+M_{x^i}(t), \quad 0<t\leq T.
\de
Thus, by the H\"older inequality and the It\^o isometry, it holds that for any $(\mathscr{F}_t)_{t\geq0}-$stopping time $\tau\leq T$ and any $\delta>0$,
\ce
\mE[|X^\e_{\tau+\delta}-X^\e_\tau|^2]&=&\sum_{i=1}^n\mE[|(X^\e_{\tau+\delta})^i-(X^\e_\tau)^i|^2]\\
&\leq&2\delta\sum_{i=1}^n\mE\left[\int_\tau^{\tau+\delta}|b^i_1(X^\e_r,Z^\e_r)|^2\dif r\right]+2\sum_{i=1}^n\mE[|M_{x^i}(\tau+\delta)-M_{x^i}(\tau)|^2]\\
&=&2\delta\sum_{i=1}^n\mE\left[\int_\tau^{\tau+\delta}|b^i_1(X^\e_r,Z^\e_r)|^2\dif r\right]+2\sum_{i=1}^n\mE\left[\int_\tau^{\tau+\delta}(\sigma_1\sigma_1^T)^{ii}(X^\e_r,Z^\e_r)\dif r\right]\\
&&+2\sum_{i=1}^n\mE\left[\int_\tau^{\tau+\delta}\int_{\mU_1}|f^i_1(X^\e_r,u)|^2\nu_1(\dif u)\dif r\right]\\
&=&2\delta\mE\left[\int_\tau^{\tau+\delta}|b_1(X^\e_r,Z^\e_r)|^2\dif r\right]+2\mE\left[\int_\tau^{\tau+\delta}\|\sigma_1(X^\e_r,Z^\e_r)\|^2\dif r\right]\\
&&+2\mE\left[\int_\tau^{\tau+\delta}\int_{\mU_1}|f_1(X^\e_r,u)|^2\nu_1(\dif u)\dif r\right]\\
&\leq&C\delta,
\de
where the last inequality is based on the condition ($\mathbf{H}^2_{b_1, \sigma_1, f_1}$), and the constant $C$ is independent of $\e$.
So,
\ce
\limsup\limits_{\delta\downarrow0}\limsup\limits_{\e\downarrow0}\sup_{\tau\leq T}\mE[|X^\e_{\tau+\delta}-X^\e_\tau|^2]=0.
\de
By the similar deduction to above, one could furthermore get
$$
\sup\limits_{\e}\sup\limits_{t\leq T}\mE[|X^\e_t|^2]<\infty.
$$
Thus, it follows from Theorem \ref{tigth} that $\{X_t^\e, t\in[0,T]\}$ is relatively weakly compact in $D([0,T],\mR^n)$.

{\bf Step 2.} We prove that the weak limit of $\{X_t^\e, t\in[0,T]\}$ is $\{X^0_t, t\in[0,T]\}$.

Taking $H(x,z)=g(x)$ in (\ref{xzmartprob}), where $g$ is a smooth and bounded function, we have that
\ce
g(X^\e_t)-g(X^\e_s)-\int_s^t(\cL^{X^\e}g)(X^\e_r,Z^\e_r)\dif r=M_g(t)-M_g(s),
\de
and then
\ce
g(X^\e_t)-g(X^\e_s)-\int_s^t(\bar{\cL}g)(X^\e_r)\dif r=\int_s^t\left[(\cL^{X^\e}g)(X^\e_r,Z^\e_r)-(\bar{\cL}g)(X^\e_r)\right]\dif r
+M_g(t)-M_g(s).
\de
Moreover, multiplying a bounded $\mathscr{F}_s$-measurable functional $\chi_s$ of the process $\{X^0_t, t\in[0,T]\}$
and taking expectation on the two hand sides of the above equality, we know
\be
\mE\left[\chi_s\left(g(X^\e_t)-g(X^\e_s)-\int_s^t(\bar{\cL}g)(X^\e_r)\dif r\right)\right]
=\mE\left[\chi_s\int_s^t\left[(\cL^{X^\e}g)(X^\e_r,Z^\e_r)-(\bar{\cL}g)(X^\e_r)\right]\dif r\right].
\label{app}
\ee

Next we compute $\lim\limits_{\e\downarrow0}\mE\left[\chi_s\int_s^t\left[(\cL^{X^\e}g)(X^\e_r,Z^\e_r)-(\bar{\cL}g)(X^\e_r)\right]\dif r\right]$.
On one hand, set
 \ce
 &\Psi(x, z, A):=\int_0^\infty[p(x;z,t,A)-\bar{p}(x;A)]\dif t,\\
 &\Psi_g(x, z):=\int_{\mR^m}\left[(\cL^{X^\e}g)(x,z')-(\bar{\cL}g)(x)\right]\Psi(x, z, \dif z'),
 \de
and then
\ce
\Psi_g(x, z)&=&\int_{\mR^m}\left[(\cL^{X^\e}g)(x,z')-(\bar{\cL}g)(x)\right]\int_0^\infty[p(x;z,t,\dif z')-\bar{p}(x;\dif z')]\dif t\\
&=&\int_0^\infty\left(\int_{\mR^m}\left[(\cL^{X^\e}g)(x,z')-(\bar{\cL}g)(x)\right][p(x;z,t,\dif z')-\bar{p}(x;\dif z')]\right)\dif t\\
&=&\int_0^\infty T^x_t[(\cL^{X^\e}g)-(\bar{\cL}g)](x,z)\dif t.
\de
Furthermore, it holds that
\be
\e(\cL^{Z^\e}\Psi_g)(x,z)&=&\int_0^\infty (\e\cL^{Z^\e}T^x_t)[(\cL^{X^\e}g)-(\bar{\cL}g)](x,z)\dif t\no\\
&=&\int_0^\infty\frac{\dif T^x_t[(\cL^{X^\e}g)-(\bar{\cL}g)](x,z)}{\dif t}\dif t\no\\
&=&\lim\limits_{t\rightarrow\infty}T^x_t[(\cL^{X^\e}g)-(\bar{\cL}g)](x,z)-[(\cL^{X^\e}g)(x,z)-(\bar{\cL}g)(x)]\no\\
&=&\lim\limits_{t\rightarrow\infty}\int_{\mR^m}\left[(\cL^{X^\e}g)(x,z')-(\bar{\cL}g)(x)\right][p(x;z,t,\dif z')-\bar{p}(x;\dif z')]\no\\
&&-[(\cL^{X^\e}g)(x,z)-(\bar{\cL}g)(x)]\no\\
&=&-[(\cL^{X^\e}g)(x,z)-(\bar{\cL}g)(x)],
\label{opeq}
\ee
where the last equality is based on {\bf Assumption 2}. On the other hand, taking $H(x,z)=\e\Psi_g(x,z)$ again in (\ref{xzmartprob}), we get
that
\ce
&&\e\Psi_g(X^\e_t,Z^\e_t)-\e\Psi_g(X^\e_s,Z^\e_s)-\e\int_s^t(\cL^{X^\e}\Psi_g)(X^\e_r,Z^\e_r)\dif r\\
&=&\int_s^t\e(\cL^{Z^\e}\Psi_g)(X^\e_r,Z^\e_r)\dif r+M_{\e\Psi_g}(t)-M_{\e\Psi_g}(s).
\de
So, by multiplying $\chi_s$ and taking expectation on the two hand sides of the above equality, it holds that
\ce
&&\e\mE\left[\chi_s\left(\Psi_g(X^\e_t,Z^\e_t)-\Psi_g(X^\e_s,Z^\e_s)-\int_s^t(\cL^{X^\e}\Psi_g)(X^\e_r,Z^\e_r)\dif r\right)\right]\\
&=&\mE\left[\chi_s\int_s^t\e(\cL^{Z^\e}\Psi_g)(X^\e_r,Z^\e_r)\dif r\right]\\
&=&-\mE\left[\chi_s\int_s^t\left[(\cL^{X^\e}g)(X^\e_r,Z^\e_r)-(\bar{\cL}g)(X^\e_r)\right]\dif r\right],
\de
where the last equality is based on (\ref{opeq}). As $\e\rightarrow0$, it is easy to see that
$$
\lim\limits_{\e\downarrow0}\mE\left[\chi_s\int_s^t\left[(\cL^{X^\e}g)(X^\e_r,Z^\e_r)-(\bar{\cL}g)(X^\e_r)\right]\dif r\right]=0,
$$
which, together with (\ref{app}), yields that
$$
\lim\limits_{\e\downarrow0}\mE\left[\chi_s\left(g(X^\e_t)-g(X^\e_s)-\int_s^t(\bar{\cL}g)(X^\e_r)\dif r\right)\right]=0.
$$
That is, the weak limit of $\{X_t^\e, t\in[0,T]\}$ is a solution of the martingale problem for ($\bar{\cL}$, $\delta_{x_0}$). Since the martingale problem for ($\bar{\cL}$, $\delta_{x_0}$) is well-posed, the weak limit of $\{X_t^\e, t\in[0,T]\}$ is $\{X^0_t, t\in[0,T]\}$.
\end{proof}

\section{Nonlinear filtering problems}\label{filpro}

In the section, we study nonlinear filtering problems for $X^\e$ and $X^0$.

For the observation process $Y^{\e}$ defined in (\ref{ypde}), i.e.
\ce
Y_t^{\e}=B_t+\int_0^th(X_s^{\e},Z^\e_s)\dif s+\int_0^t\int_{\mU_3}f_3(s,u)\tilde{N}_{\lambda}(\dif s, \dif u)+\int_0^t\int_{\mU\setminus\mU_3}g_3(s,u)N_{\lambda}(\dif s, \dif u),
\de
we assume:

\begin{enumerate}[\bf{Assumption 4.}]
\item
\end{enumerate}
\begin{enumerate}[\bf{(i)}]
\item $h$ is bounded and
$$
\int_0^T\int_{\mU_3}|f_3(s,u)|^2\nu_3(\dif u)\dif s<\infty.
$$
\end{enumerate}
\begin{enumerate}[\bf{(ii)}]
\item There exists a positive function $\hat{L}(u)$ satisfying
\ce
\int_{\mU_3}\frac{\left(1-\hat{L}(u)\right)^2}{\hat{L}(u)}\nu_3(\dif u)<\infty
\de
such that $0<\hat{l}\leq \hat{L}(u)<\lambda(t,x,u)<1$ for $u\in\mU_3$, where $\hat{l}$ is a constant.
\end{enumerate}

Under {\bf Assumption 4.(i)}, $Y^{\e}$ is well defined. And set
\ce
(\Lambda^\e_t)^{-1}:&=&\exp\bigg\{-\int_0^th(X_s^{\e},Z^\e_s)^i\dif B^i_s-\frac{1}{2}\int_0^t
\left|h(X_s^{\e},Z^\e_s)\right|^2\dif s
-\int_0^t\int_{\mU_3}\log\lambda(s,X^{\e}_{s-},u)N_{\lambda}(\dif s, \dif u)\\
&&\quad\qquad -\int_0^t\int_{\mU_3}(1-\lambda(s,X^\e_s,u))\nu_3(\dif u)\dif s\bigg\}.
\de
So, under {\bf Assumption 4.(ii)}, it holds that
\ce
\mE\left[\exp\left\{\int_0^T\int_{\mU_3}\frac{\left(1-\lambda(s,X^\e_s,u)\right)^2}{\lambda(s,X^\e_s,u)}
\nu_3(\dif u)\dif s\right\}\right]
\leq\exp\left\{\int_0^T\int_{\mU_3}\frac{\left(1-\hat{L}(u)\right)^2}{\hat{L}(u)}\nu_3(\dif u)\dif s\right\}
<\infty.
\de
Thus, by the same deduction as that in \cite{qd}, we know that $(\Lambda^\e_t)^{-1}$
is an exponential martingale. By use of $(\Lambda^\e_t)^{-1}$, one could define a measure $\mP^\e$ via
$$
\frac{\dif \mP^\e}{\dif \mP}=(\Lambda^\e_T)^{-1}.
$$
By the Girsanov theorem for Brownian motions and random measures, we can obtain that under the measure
$\mP^\e$, $\bar{B}_t:=B_t+\int_0^th(X_s^{\e},Z^\e_s)\dif s$ is a Brownian motion and $N_{\lambda}((0,t],\dif u)$
is a Poisson random measure with the predictable compensator $t\nu_3(\dif u)$.

Next, we rewrite $\Lambda^\e_t$ as
\ce
\Lambda^\e_t&=&\exp\bigg\{\int_0^th(X_s^{\e},Z^\e_s)^i\dif \bar{B}^i_s-\frac{1}{2}\int_0^t
\left|h(X_s^{\e},Z^\e_s)\right|^2\dif s
+\int_0^t\int_{\mU_3}\log\lambda(s,X^{\e}_{s-},u)N_{\lambda}(\dif s, \dif u)\\
&&\quad\qquad +\int_0^t\int_{\mU_3}(1-\lambda(s,X^\e_s,u))\nu_3(\dif u)\dif s\bigg\}.
\de
Define
\ce
&&\rho^{\e}_t(\psi):=\mE^{\mP^\e}[\psi(X^{\e}_t)\Lambda^\e_t|\mathscr{F}_t^{Y^{\e}}], \\
&&\pi^{\e}_t(\psi):=\mE[\psi(X^{\e}_t)|\mathscr{F}_t^{Y^{\e}}], \qquad \psi\in\cB(\mR^n),
\de
where $\mE^{\mP^\e}$ denotes the expectation under the measure $\mP^\e$ and $\mathscr{F}_t^{Y^{\e}}$ stands for the usual augmentation of the $\sigma$-algebra generated by $\{Y^{\e}_s, 0\leq s\leq t\}$. $ \rho_t^\e$ and $\pi^{\e}_t$ are called the nonnormalized filtering and the normalized filtering of $X_t^\e$ with respect to $\mathscr{F}_t^{Y^{\e}}$, respectively. And then by the Kallianpur-Striebel formula it holds that
\ce
\pi^{\e}_t(\psi)=\frac{\rho^{\e}_t(\psi)}{\rho^{\e}_t(1)}.
\de

Set
\ce
\bar{h}(x):=\int_{\mR^m}h(x,z)\bar{p}(x,\dif z),
\de
and then $\bar{h}$ is an averaged version of $h$. So, we make use of $\bar{h}$ to define
\ce
\bar{\Lambda}_t&:=&\exp\bigg\{\int_0^t\bar{h}(X_s^0)^i\dif \bar{B}^i_s-\frac{1}{2}\int_0^t
\left|\bar{h}(X_s^0)\right|^2\dif s
+\int_0^t\int_{\mU_3}\log\lambda(s,X^0_{s-},u)N_{\lambda}(\dif s, \dif u)\\
&&\quad\qquad +\int_0^t\int_{\mU_3}(1-\lambda(s,X^0_s,u))\nu_3(\dif u)\dif s\bigg\},
\de
and furthermore
\ce
\rho^0_t(\psi)&:=&\mE^{\mP^\e}[\psi(X^0_t)\bar{\Lambda}_t|\mathscr{F}_t^{Y^{\e}}],\\
\pi^0_t(\psi)&:=&\frac{\rho^0_t(\psi)}{\rho^0_t(1)}.
\de
And then we study the relation between $\pi^0_t$ and $\pi^{\e}_t$ as $\e\rightarrow0$ in the next section.

At the first look, it is more reasonable to define the limit observable process
\ce
Y_t^{0}:=\int_0^t\bar{h}(X_s^0)\dif s+B_t+\int_0^t\int_{\mU_3}f_3(s,u)\tilde{\bar{N}}_{\lambda}(\dif s, \dif u)+\int_0^t\int_{\mU\setminus\mU_3}g_3(s,u)\bar{N}_{\lambda}(\dif s, \dif u),
\de
where $\bar{N}_{\lambda}(\dif t, \dif u)$ is a random measure with a predictable compensator \\$\lambda(t,X_t^0,u)\dif t\nu_3(\dif u)$,
and the corresponding nonlinear filtering
\ce
\mP^0_t(\psi):=\mE[\psi(X^0_t)|\mathscr{F}_t^{Y^0}],
\de
and discuss the relation between $\mP^0_t$ and $\pi^{\e}_t$ as $\e\rightarrow0$. In fact, since
$X_t^0$ couldn't be obtained genuinely, $Y_t^{0}$ is not observable. However, should such homogenized observation be
available, using it would lead to loss of information for estimating the signal
compared to using the actual observation. Therefore, we only consider $X_t^0$ under
$\mathscr{F}_t^{Y^{\e}}$.

\section{Convergence of nonlinear filterings}\label{confil}

In the section, we prove that $\pi^{\e}_t$ converges weakly to $\pi^0_t$ as $\e\rightarrow0$ for any $t\in[0,T]$. 
Firstly, let us prove two key lemmas.

\bl\label{rho1}
Suppose that $h, \lambda$ satisfy {\bf Assumption 4.}. Then it holds that for any $t\in[0,T]$,
\be
(\rho^0_t(1))^{-1}<\infty, \qquad \mP\ a.s..
\label{roes}
\ee
\el
\begin{proof}
By the H\"older inequality, it holds that
$$
\mE(\rho^0_t(1))^{-1}=\mE^{\mP^\e}\left(\rho^0_t(1)\right)^{-1}\Lambda^\e_T\leq\left(\mE^{\mP^\e}(\rho^0_t(1))^{-2}\right)^{1/2}
\left(\mE^{\mP^\e}(\Lambda^\e_T)^2\right)^{1/2}.
$$
Let us firstly estimate $\mE^{\mP^\e}(\rho^0_t(1))^{-2}$. Note that $\rho^0_t(1)=\mE^{\mP^\e}[\bar{\Lambda}_t|\mathscr{F}_t^{Y^{\e}}]$
and $x\mapsto x^{-2}$ is convex. Thus, we know by the Jensen inequality that
$$
\mE^{\mP^\e}(\rho^0_t(1))^{-2}=\mE^{\mP^\e}(\mE^{\mP^\e}[\bar{\Lambda}_t|\mathscr{F}_t^{Y^{\e}}])^{-2}
\leq\mE^{\mP^\e}[\mE^{\mP^\e}[(\bar{\Lambda}_t)^{-2}|\mathscr{F}_t^{Y^{\e}}]]=\mE^{\mP^\e}(\bar{\Lambda}_t)^{-2}.
$$
So, we estimate $\mE^{\mP^\e}(\bar{\Lambda}_t)^{-2}$. Applying the It\^o formula to $(\bar{\Lambda}_t)^{-1}$, one could obtain that
\ce
(\bar{\Lambda}_t)^{-1}&=&1+\int_0^t(\bar{\Lambda}_s)^{-1}|\bar{h}(X_s^0)|^2\dif s+\int_0^t\int_{\mU_3}(\bar{\Lambda}_s)^{-1}\frac{(1-\lambda(s,X_s^0,u))^2}{\lambda(s,X_s^0,u)}\nu_3(\dif u)\dif s\\
&&-\int_0^t(\bar{\Lambda}_s)^{-1}\bar{h}(X_s^0)^i\dif \bar{B}^i_s+\int_0^t\int_{\mU_3}(\bar{\Lambda}_s)^{-1}\frac{1-\lambda(s,X_s^0,u)}{\lambda(s,X_s^0,u)}\tilde{N}_{\lambda}(\dif s, \dif u).
\de
Furthermore, it follows from the H\"older inequality and the It\^o isometry that
\ce
\mE^{\mP^\e}(\bar{\Lambda}_t)^{-2}&\leq&5+5\mE^{\mP^\e}\left|\int_0^t(\bar{\Lambda}_s)^{-1}|\bar{h}(X_s^0)|^2\dif s\right|^2+5\mE^{\mP^\e}\left|\int_0^t(\bar{\Lambda}_s)^{-1}\bar{h}(X_s^0)^i\dif \bar{B}^i_s\right|^2\\
&&+5\mE^{\mP^\e}\left|\int_0^t\int_{\mU_3}(\bar{\Lambda}_s)^{-1}\frac{(1-\lambda(s,X_s^0,u))^2}{\lambda(s,X_s^0,u)}\nu_3(\dif u)\dif s\right|^2\\
&&+5\mE^{\mP^\e}\left|\int_0^t\int_{\mU_3}(\bar{\Lambda}_s)^{-1}\frac{1-\lambda(s,X_s^0,u)}{\lambda(s,X_s^0,u)}\tilde{N}_{\lambda}(\dif s, \dif u)\right|^2\\
&\leq&5+5T\mE^{\mP^\e}\int_0^t(\bar{\Lambda}_s)^{-2}|\bar{h}(X_s^0)|^4\dif s+5\mE^{\mP^\e}\int_0^t(\bar{\Lambda}_s)^{-2}|\bar{h}(X_s^0)|^2\dif s\\
&&+5T\mE^{\mP^\e}\int_0^t(\bar{\Lambda}_s)^{-2}\left|\int_{\mU_3}\frac{(1-\lambda(s,X_s^0,u))^2}{\lambda(s,X_s^0,u)}\nu_3(\dif u)\right|^2\dif s\\
&&+5\mE^{\mP^\e}\int_0^t\int_{\mU_3}(\bar{\Lambda}_s)^{-2}\frac{(1-\lambda(s,X_s^0,u))^2}{\lambda(s,X_s^0,u)^2}\nu_3(\dif u)\dif s\\
&\leq&5+C\int_0^t\mE^{\mP^\e}(\bar{\Lambda}_s)^{-2}\dif s,
\de
where the last step is based on {\bf Assumption 4.}. The Gronwall inequality admits us to have $\mE^{\mP^\e}(\bar{\Lambda}_t)^{-2}<\infty$.

Next, we deal with $\mE^{\mP^\e}(\Lambda^\e_T)^2$. Applying the It\^o formula to $\Lambda^\e_t$, one can obtain that
\be
\Lambda^\e_t=1+\int_0^t\Lambda^\e_sh(X^\e_s,Z^\e_s)^i\dif \bar{B}^i_s+\int_0^t\int_{\mU_3}\Lambda^\e_{s-}(\lambda(s,X^\e_{s-},u)-1)\tilde{N}_{\lambda}(\dif s, \dif u).
\label{lme}
\ee
Thus, by the similar deduction to $\mE^{\mP^\e}(\bar{\Lambda}_t)^{-2}$ it holds that $\mE^{\mP^\e}(\Lambda^\e_T)^2<\infty$.

In conclusion, $\mE(\rho^0_t(1))^{-1}<\infty$ and then (\ref{roes}) is right. The proof is completed.
\end{proof}

\bl\label{tigh}
Under {\bf Assumption 4.}, $\{\rho^{\e}_t, t\in[0,T]\}$ is relatively weakly compact in $D([0,T], \cM(\mR^n))$.
\el
\begin{proof}
First of all, we explain $\rho^{\e}_t\in\cM(\mR^n)$ for $\omega\in\Omega, t\in[0,T]$. Note that $\rho^{\e}_t(\mR^n)=\rho^{\e}_t(1_{\mR^n})=\rho^{\e}_t(1)$.
And then by the H\"older inequality, it holds that
$$
\mE\rho^{\e}_t(\mR^n)=\mE\rho^{\e}_t(1)=\mE^{\mP^\e}[\rho^{\e}_t(1)\Lambda^\e_T]
\leq\left(\mE^{\mP^\e}(\rho^{\e}_t(1))^2\right)^{1/2}\left(\mE^{\mP^\e}(\Lambda^\e_T)^2\right)^{1/2}.
$$
On one hand, the Jensen inequality admits us to obtain that
$$
\mE^{\mP^\e}(\rho^{\e}_t(1))^2=\mE^{\mP^\e}[\mE^{\mP^\e}[\Lambda^\e_t|\mathscr{F}_t^{Y^{\e}}]]^2
\leq\mE^{\mP^\e}[\mE^{\mP^\e}[(\Lambda^\e_t)^2|\mathscr{F}_t^{Y^{\e}}]]=\mE^{\mP^\e}(\Lambda^\e_t)^2.
$$
By the proof of Lemma \ref{rho1}, one could get $\mE^{\mP^\e}(\Lambda^\e_t)^2<\infty$ and then $\mE^{\mP^\e}(\rho^{\e}_t(1))^2<\infty$. On
the other hand, it follows from the proof of Lemma \ref{rho1} that $\mE^{\mP^\e}(\Lambda^\e_T)^2<\infty$. Thus, $\mE\rho^{\e}_t(\mR^n)<\infty$ and then $\rho^{\e}_t(\mR^n)<\infty$ a.s. $\mP$. In addition, it is easy to justify other measure properties of $\rho^{\e}_t$ by means of properties of conditional expectations.

Next, we deduce the equation for $\rho^{\e}_t$. For $\psi\in \cC^2_b(\mR^n)$, applying the It\^o formula to $\psi(X^\e_t)$, one can have that
\ce
\psi(X^\e_t)&=&\psi(X^\e_0)+\int_0^t(\cL^{X^\e}\psi)(X^\e_s,Z^\e_s)\dif s+\int_0^t(\nabla\psi)(X^\e_s)\sigma_1(X^\e_s,Z^\e_s)\dif V_s\\
&&+\int_0^t\int_{\mU_1}[\psi(X^\e_{s-}+f_1(X^\e_{s-}, u))-\psi(X^\e_{s-})]\tilde{N}_{p_1}(\dif s, \dif u).
\de
Note that $\Lambda^\e_t$ satisfies (\ref{lme}). So, it follows from the It\^o formula that
\ce
\psi(X^\e_t)\Lambda^\e_t&=&\psi(X^\e_0)+\int_0^t\psi(X^\e_s)\Lambda^\e_sh(X_s^{\e},Z^\e_s)^i\dif \bar{B}^i_s\\
&&+\int_0^t\int_{\mU_3}\psi(X^\e_{s-})\Lambda^\e_{s-}(\lambda(s,X^\e_{s-},u)-1)\tilde{N}_{\lambda}(\dif s, \dif u)\\
&&+\int_0^t\Lambda^\e_s(\cL^{X^\e}\psi)(X^\e_s,Z^\e_s)\dif s+\int_0^t\Lambda^\e_s(\nabla\psi)(X^\e_s)\sigma_1(X^\e_s,Z^\e_s)\dif V_s\\
&&+\int_0^t\int_{\mU_1}\Lambda^\e_{s-}[\psi(X^\e_{s-}+f_1(X^\e_{s-}, u))-\psi(X^\e_{s-})]\tilde{N}_{p_1}(\dif s, \dif u).
\de
Taking the conditional expectation with respect to $\mathscr{F}_t^{Y^{\e}}$ under $\mP^\e$ on two hand sides of the above
equality, one could obtain  that
\ce
\mE^{\mP^\e}[\psi(X^\e_t)\Lambda^\e_t|\mathscr{F}_t^{Y^{\e}}]&=&\mE^{\mP^\e}[\psi(X^\e_0)|\mathscr{F}_0^{Y^{\e}}]
+\int_0^t\mE^{\mP^\e}[\psi(X^\e_s)\Lambda^\e_sh(X_s^{\e},Z^\e_s)^i|\mathscr{F}_s^{Y^{\e}}]\dif \bar{B}^i_s\\
&&+\int_0^t\int_{\mU_3}\mE^{\mP^\e}[\psi(X^\e_{s-})\Lambda^\e_{s-}(\lambda(s,X^\e_{s-},u)-1)|\mathscr{F}_s^{Y^{\e}}]\tilde{N}_{\lambda}(\dif s, \dif u)\\
&&+\int_0^t\mE^{\mP^\e}[\Lambda^\e_s(\cL^{X^\e}\psi)(X^\e_s,Z^\e_s)|\mathscr{F}_s^{Y^{\e}}]\dif s,
\de
i.e.
\be
\rho^{\e}_t(\psi)&=&\rho^{\e}_0(\psi)+\int_0^t\rho^{\e}_s\(\big(\cL^{X^\e}\psi\big)(\cdot,Z^\e_s)\)\dif s+\int_0^t\rho^{\e}_s\(\psi h(\cdot,Z^\e_s)^i\)\dif \bar{B}^i_s\no\\
&&+\int_0^t\int_{\mU_3}\rho^{\e}_s\(\psi(\lambda(s,\cdot,u)-1)\)\tilde{N}_{\lambda}(\dif s, \dif u).
\label{zakai1}
\ee
For the detailed deduction of the above equation, please refer to the proof of Theorem 3.1 in \cite{qd}.

For any $(\mathscr{F}_t)_{t\geq0}-$stopping time $\tau\leq T$ and any $\delta>0$, we compute $\mE|\rho^{\e}_{\tau+\delta}(\psi)-\rho^{\e}_{\tau}(\psi)|$.
It follows from the H\"older inequality that
\ce
\mE|\rho^{\e}_{\tau+\delta}(\psi)-\rho^{\e}_{\tau}(\psi)|=\mE^{\mP^\e}|\rho^{\e}_{\tau+\delta}(\psi)-\rho^{\e}_{\tau}(\psi)|\Lambda^\e_T
\leq\left(\mE^{\mP^\e}|\rho^{\e}_{\tau+\delta}(\psi)-\rho^{\e}_{\tau}(\psi)|^2\right)^{1/2}\left(\mE^{\mP^\e}(\Lambda^\e_T)^2\right)^{1/2}.
\de
Since $\mE^{\mP^\e}(\Lambda^\e_T)^2<C$, which has been proved in Lemma \ref{rho1}, we only consider $\mE^{\mP^\e}|\rho^{\e}_{\tau+\delta}(\psi)-\rho^{\e}_{\tau}(\psi)|^2$. The H\"older inequality and
the It\^o isometry admit us to get
\ce
\mE^{\mP^\e}|\rho^{\e}_{\tau+\delta}(\psi)-\rho^{\e}_{\tau}(\psi)|^2&\leq&3\mE^{\mP^\e}\left|\int_{\tau}^{\tau+\delta}\rho^{\e}_s\(\big(\cL^{X^\e}\psi\big)(\cdot,Z^\e_s)\)\dif s\right|^2\\
&&+3\mE^{\mP^\e}\left|\int_{\tau}^{\tau+\delta}\rho^{\e}_s\(\psi h(\cdot,Z^\e_s)^i\)\dif \bar{B}^i_s\right|^2\\
&&+3\mE^{\mP^\e}\left|\int_{\tau}^{\tau+\delta}\int_{\mU_3}\rho^{\e}_s\(\psi(\lambda(s,\cdot,u)-1)\)\tilde{N}_{\lambda}(\dif s, \dif u)\right|^2\\
&\leq&3\delta\mE^{\mP^\e}\int_{\tau}^{\tau+\delta}\left|\rho^{\e}_s\(\big(\cL^{X^\e}\psi\big)(\cdot,Z^\e_s)\)\right|^2\dif s\\
&&+3\mE^{\mP^\e}\int_{\tau}^{\tau+\delta}\left|\rho^{\e}_s\(\psi h(\cdot,Z^\e_s)^i\)\right|^2\dif s\\
&&+3\mE^{\mP^\e}\int_{\tau}^{\tau+\delta}\int_{\mU_3}\left|\rho^{\e}_s\(\psi(\lambda(s,\cdot,u)-1)\)\right|^2\nu_3(\dif u)\dif s\\
&=:&I_1+I_2+I_3.
\de
For $I_1$, by the Jensen inequality and ($\mathbf{H}^2_{b_1, \sigma_1, f_1}$), it holds that
\ce
I_1&=&3\delta\mE^{\mP^\e}\int_{\tau}^{\tau+\delta}\left|\mE^{\mP^\e}[\Lambda^\e_s(\cL^{X^\e}\psi)(X^\e_s,Z^\e_s)|\mathscr{F}_s^{Y^{\e}}]\right|^2\dif s\\
&\leq&3\delta\mE^{\mP^\e}\int_{\tau}^{\tau+\delta}\mE^{\mP^\e}[(\Lambda^\e_s)^2\left|(\cL^{X^\e}\psi)(X^\e_s,Z^\e_s)\right|^2|\mathscr{F}_s^{Y^{\e}}]\dif s\\
&\leq&3C\delta\mE^{\mP^\e}\int_{\tau}^{\tau+\delta}\mE^{\mP^\e}[(\Lambda^\e_s)^2|\mathscr{F}_s^{Y^{\e}}]\dif s\\
&=&3C\delta\mE^{\mP^\e}\int_{0}^{\delta}\mE^{\mP^\e}[(\Lambda^\e_{\tau+s})^2|\mathscr{F}_{\tau+s}^{Y^{\e}}]\dif s\\
&=&3C\delta\int_{0}^{\delta}\mE^{\mP^\e}[\mE^{\mP^\e}[(\Lambda^\e_{\tau+s})^2|\mathscr{F}_{\tau+s}^{Y^{\e}}]]\dif s\\
&=&3C\delta\int_{0}^{\delta}\mE^{\mP^\e}[(\Lambda^\e_{\tau+s})^2]\dif s\\
&\leq&3C\delta^2,
\de
where $C$ is independent of $\e, \delta$. By the same deduction as $I_1$, we get that $I_2+I_3\leq C\delta$. Thus,
\be
\limsup_{\delta\downarrow0}\limsup_{\e\downarrow0}\sup\limits_{\tau\leq T}\mE|\rho^{\e}_{\tau+\delta}(\psi)-\rho^{\e}_{\tau}(\psi)|=0.
\label{eqco}
\ee

Based on the similar calculation to above, it holds that
\be
\sup\limits_{\e}\sup\limits_{t\leq T}\mE|\rho^{\e}_t(\psi)|<\infty.
\label{unbo}
\ee
So, combining (\ref{unbo}) with (\ref{eqco}), we know from Theorem 5.1 in \cite{kus} that $\{\rho^{\e}_t(\psi), t\in[0,T]\}$ is
relatively weakly compact in $D([0,T], \mR)$. Moreover, Theorem 6.2 in \cite{kus} admits us to obtain that
$\{\rho^{\e}_t, t\in[0,T]\}$ is relatively weakly compact in $D([0,T], \cM(\mR^n))$.
\end{proof}

To attain the convergence of $\pi^{\e}_t$ to $\pi^0_t$  as $\e\rightarrow0$, we assume more:

\medspace

{\bf Assumption 5.}
\begin{enumerate}[]  
\item $\{Z_{\e t}^{\e}, t\in[0,T]\}$ is tight.
\end{enumerate}

Now, it is the position to state the main result in the section.

\bt\label{filcon}
Under {\bf Assumption 1.-5.}, $\pi^{\e}_t$ converges weakly to $\pi^0_t$ as $\e\rightarrow0$ for any $t\in[0,T]$.
\et
\begin{proof}
For $\psi\in\cC^2_b(\mR^n)$, it holds that
\ce
\pi^{\e}_t(\psi)-\pi^0_t(\psi)=\frac{\rho^{\e}_t(\psi)-\rho^0_t(\psi)}{\rho^0_t(1)}-\pi^{\e}_t(\psi)\frac{\rho^{\e}_t(1)-\rho^0_t(1)}{\rho^0_t(1)}.
\de
Thus, in order to prove $\pi^{\e}_t(\psi)-\pi^0_t(\psi)$ converges weakly to $0$, by Lemma \ref{rho1} and the conditional expectation property of $\pi^{\e}_t(\psi)$, we only need to show that $\rho^{\e}_t(\psi)$ converges weakly to $\rho^0_t(\psi)$ as $\e\rightarrow0$. Besides, by Lemma \ref{tigh}, there exist a weakly convergence subsequence $\{\rho^{\e_k}_t, k\in\mN\}$ and a measure-valued process $\bar{\rho}_t$ such that $\rho^{\e_k}_t(\psi)$ converges weakly to $\bar{\rho}_t(\psi)$ as $k\rightarrow\infty$. Therefore, we just need to prove that for $t\in[0,T]$
\ce
\bar{\rho}_t(\psi)=\rho^0_t(\psi), \qquad a.s. \mP.
\de

Next, we are devoted to looking for the equations which $\bar{\rho}_t(\psi), \rho^0_t(\psi)$ solve. First of all, let us deduce the equation which $\bar{\rho}_t(\psi)$ satisfies. Note that $\rho^{\e}_t(\psi)$ solves
Eq.(\ref{zakai1}). And then we study the weak limits of three integrals in Eq.(\ref{zakai1}). 

For the first integral of Eq.(\ref{zakai1}), we consider the limit of $\rho^{\e_k}_s\(\big(\cL^{X^{\e_k}}\psi\big)(\cdot,Z^{\e_k}_s)\)-\bar{\rho}_s\(\bar{\cL}\psi\)$ as $k\rightarrow\infty$. By the simple calculation, it holds that
\ce
\rho^{\e_k}_s\(\big(\cL^{X^{\e_k}}\psi\big)(\cdot,Z^{\e_k}_s)\)-\bar{\rho}_s\(\bar{\cL}\psi\)
&=&\rho^{\e_k}_s\(\big(\cL^{X^{\e_k}}\psi\big)(\cdot,Z^{\e_k}_s)\)-\rho^{\e_k}_s\(\bar{\cL}\psi\)\\
&&+\rho^{\e_k}_s\(\bar{\cL}\psi\)-\bar{\rho}_s\(\bar{\cL}\psi\)\\
&=&\rho^{\e_k}_s\(\big(\cL^{X^{\e_k}}\psi\big)(\cdot,Z^{\e_k}_s)-\bar{\cL}\psi\)\\
&&+\left[\rho^{\e_k}_s\(\bar{\cL}\psi\)-\bar{\rho}_s\(\bar{\cL}\psi\)\right]\\
&=:&I_1+I_2.
\de
For $I_1$, one know that
\ce
I_1&=&\mE^{\mP^{\e_k}}\left[\Lambda^{\e_k}_s\left(\big(\cL^{X^{\e_k}}\psi\big)(X^{\e_k}_s,Z^{\e_k}_s)-(\bar{\cL}\psi)(X^{\e_k}_s)\right)|\mathscr{F}_s^{Y^{\e_k}}\right]\\
&=&\mE^{\mP^{\e_k}}\left[\Lambda^{\e_k}_s\frac{\partial\psi}{\partial x_i}(X^{\e_k}_s)\left[b_1^i(X^{\e_k}_s,Z^{\e_k}_s)-\bar{b}^i_1(X^{\e_k}_s)\right]|\mathscr{F}_s^{Y^{\e_k}}\right]\\
&&+\frac{1}{2}\mE^{\mP^{\e_k}}\left[\Lambda^{\e_k}_s\frac{\partial^2\psi}{\partial x_i\partial x_j}(X^{\e_k}_s)\left[(\sigma_1\sigma_1^T)^{ij}(X^{\e_k}_s,Z^{\e_k}_s)-(\bar{\sigma}_1\bar{\sigma}_1^T)^{ij}(X^{\e_k}_s)\right]|\mathscr{F}_s^{Y^{\e_k}}\right]\\
&=:&I_{11}+I_{12}.
\de
Let us deal with $I_{11}$. Since 
\ce
&&\lim\limits_{n\rightarrow\infty}\sum_{j=0}^{n-1}\mE^{\mP^{\e_k}}\left[\Lambda^{\e_k}_{(j+1)t/n}\frac{\partial\psi}{\partial x_i}(X^{\e_k}_{(j+1)t/n})\left[b_1^i(X^{\e_k}_{(j+1)t/n},Z^{\e_k}_s)-\bar{b}^i_1(X^{\e_k}_{(j+1)t/n})\right]\bigg|\mathscr{F}_s^{Y^{\e_k}}\right]I_{(jt/n,(j+1)t/n]}(s)\\
&=&\mE^{\mP^{\e_k}}\left[\Lambda^{\e_k}_s\frac{\partial\psi}{\partial x_i}(X^{\e_k}_s)\left[b_1^i(X^{\e_k}_s,Z^{\e_k}_s)-\bar{b}^i_1(X^{\e_k}_s)\right]\bigg|\mathscr{F}_s^{Y^{\e_k}}\right], a.s.\mP,
\de
we only consider $\mE^{\mP^{\e_k}}\left[\Lambda^{\e_k}_{(j+1)t/n}\frac{\partial\psi}{\partial x_i}(X^{\e_k}_{(j+1)t/n})\left[b_1^i(X^{\e_k}_{(j+1)t/n},Z^{\e_k}_s)-\bar{b}^i_1(X^{\e_k}_{(j+1)t/n})\right]\bigg|\mathscr{F}_s^{Y^{\e_k}}\right]$ for $s\in(jt/n,(j+1)t/n]$. Based on independence of $X_{\cdot}^{\e_k}, Z_{\cdot}^{\e_k}$ and $Y_{\cdot}^{\e_k}$ under $\mP^{\e_k}$, it holds that
\ce
&&\mE^{\mP^{\e_k}}\left[\Lambda^{\e_k}_{(j+1)t/n}\frac{\partial\psi}{\partial x_i}(X^{\e_k}_{(j+1)t/n})\left[b_1^i(X^{\e_k}_{(j+1)t/n},Z^{\e_k}_s)-\bar{b}^i_1(X^{\e_k}_{(j+1)t/n})\right]\bigg|\mathscr{F}_s^{Y^{\e_k}}\right]\\
&=&\mE^{\mP^{\e_k}}\Bigg[\Lambda^{\e_k}_{(j+1)t/n}\mE^{\mP^{\e_k}}\Big[\frac{\partial\psi}{\partial x_i}(X^{\e_k}_{(j+1)t/n})\left[b_1^i(X^{\e_k}_{(j+1)t/n},Z^{\e_k}_s)-\bar{b}^i_1(X^{\e_k}_{(j+1)t/n})\right]\\
&&\quad\quad\bigg|X^{\e_k}_{s+\e_k\iota},Z^{\e_k}_{s-\e_k\iota}\Big]\bigg|\mathscr{F}_s^{Y^{\e_k}}\Bigg],
\de
where $\iota$ is a positive integer such that $s-\e_k\iota>0$ and $s+\e_k\iota\leq(j+1)t/n$. Thus, we compute 
$$
\mE^{\mP^{\e_k}}\left[\frac{\partial\psi}{\partial x_i}(X^{\e_k}_{(j+1)t/n})\left[b_1^i(X^{\e_k}_{(j+1)t/n},Z^{\e_k}_s)-\bar{b}^i_1(X^{\e_k}_{(j+1)t/n})\right]\bigg|X^{\e_k}_{s+\e_k\iota},Z^{\e_k}_{s-\e_k\iota}\right].
$$ 
On one side, it is easy to see that
\ce
&&\mE^{\mP^{\e_k}}\left[\frac{\partial\psi}{\partial x_i}(X^{\e_k}_{(j+1)t/n})\left[b_1^i(X^{\e_k}_{(j+1)t/n},Z^{\e_k}_s)-\bar{b}^i_1(X^{\e_k}_{(j+1)t/n})\right]\bigg|X^{\e_k}_{s+\e_k\iota},Z^{\e_k}_{s-\e_k\iota}\right]\\
&&-\int_{\mR^m}\frac{\partial\psi}{\partial x_i}(X^{\e_k}_{(j+1)t/n})\left[b_1^i(X^{\e_k}_{(j+1)t/n},z)-\bar{b}^i_1(X^{\e_k}_{(j+1)t/n})\right]\bar{p}(X^{\e_k}_{s},\dif z)\\
&=&\Bigg(\mE^{\mP^{\e_k}}\left[\frac{\partial\psi}{\partial x_i}(X^{\e_k}_{(j+1)t/n})\left[b_1^i(X^{\e_k}_{(j+1)t/n},Z^{\e_k}_s)-\bar{b}^i_1(X^{\e_k}_{(j+1)t/n})\right]\bigg|X^{\e_k}_{s+\e_k\iota},Z^{\e_k}_{s-\e_k\iota}\right]\\
&&-\int_{\mR^m}\frac{\partial\psi}{\partial x_i}(X^{\e_k}_{(j+1)t/n})\left[b_1^i(X^{\e_k}_{(j+1)t/n},z)-\bar{b}^i_1(X^{\e_k}_{(j+1)t/n})\right]p(X^{\e_k}_{s+\e_k\iota}; Z^{\e_k}_{s-\e_k\iota}, \iota, \dif z)\Bigg)\\
&&+\Bigg(\int_{\mR^m}\frac{\partial\psi}{\partial x_i}(X^{\e_k}_{(j+1)t/n})\left[b_1^i(X^{\e_k}_{(j+1)t/n},z)-\bar{b}^i_1(X^{\e_k}_{(j+1)t/n})\right]p(X^{\e_k}_{s+\e_k\iota}; Z^{\e_k}_{s-\e_k\iota}, \iota, \dif z)\\
&&-\int_{\mR^m}\frac{\partial\psi}{\partial x_i}(X^{\e_k}_{(j+1)t/n})\left[b_1^i(X^{\e_k}_{(j+1)t/n},z)-\bar{b}^i_1(X^{\e_k}_{(j+1)t/n})\right]\bar{p}(X^{\e_k}_{s},\dif z)\Bigg)\\
&=:&I_{111}+I_{112}.
\de
Based on tightness of $\{(X_t^{\e}, Z_{\e t}^{\e}), t\in[0,T]\}$ and ($\mathbf{H}^1_{b_1, \sigma_1, f_1}$), it holds that $\lim\limits_{k\rightarrow\infty}I_{111}=0$. And it follows from the definition of $p(X^{\e_k}_{s+\e_k\iota}; Z^{\e_k}_{s-\e_k\iota}, \iota, \dif z)$ and $\bar{p}(X^{\e_k}_{s},\dif z)$ that $\lim\limits_{\iota\rightarrow\infty}I_{112}=0$. On the other side, the dominated convergence theorem admits us to obtain that
\ce
&&\lim\limits_{n\rightarrow\infty}\int_{\mR^m}\frac{\partial\psi}{\partial x_i}(X^{\e_k}_{(j+1)t/n})\left[b_1^i(X^{\e_k}_{(j+1)t/n},z)-\bar{b}^i_1(X^{\e_k}_{(j+1)t/n})\right]\bar{p}(X^{\e_k}_{s},\dif z)\\
&=&\int_{\mR^m}\frac{\partial\psi}{\partial x_i}(X^{\e_k}_s)\left[b_1^i(X^{\e_k}_s,z)-\bar{b}^i_1(X^{\e_k}_s)\right]\bar{p}(X^{\e_k}_s,\dif z)=0.
\de
Thus, again by the dominated convergence theorem, it holds that $\lim\limits_{k\rightarrow\infty}I_{11}=0$.

By the same deduction as that for $I_{11}$, it holds that $I_{12}$ goes to zero a.s. as $k\rightarrow\infty$. Thus, $I_1$ converges to zero as $k\rightarrow\infty$,
which together with weak convergence of $I_2$ to zero as $k\rightarrow\infty$  yields that $\rho^{\e_k}_s\(\big(\cL^{X^{\e_k}}\psi\big)(\cdot,Z^{\e_k}_s)\)$ converges weakly to $\bar{\rho}_s\(\bar{\cL}\psi\)$ as $k\rightarrow\infty$. 

Besides, set 
\ce
&&\left(\rho^{\e_k}_s\(\big(\cL^{X^{\e_k}}\psi\big)(\cdot,Z^{\e_k}_s)\)\right)^{(n)}:=\sum_{j=0}^{n-1}\rho^{\e_k}_{(j+1)t/n}\(\big(\cL^{X^{\e_k}}\psi\big)(\cdot,Z^{\e_k}_{(j+1)t/n})\)I_{(jt/n, (j+1)t/n]}(s),\\
&&\left(\bar{\rho}_s\big(\bar{\cL}\psi\big)\right)^{(n)}:=\sum_{j=0}^{n-1}\bar{\rho}_{(j+1)t/n}\big(\bar{\cL}\psi\big)I_{(jt/n, (j+1)t/n]}(s),
\de
and then 
\ce
&&\lim\limits_{n\rightarrow\infty}\left(\rho^{\e_k}_s\(\big(\cL^{X^{\e_k}}\psi\big)(\cdot,Z^{\e_k}_s)\)\right)^{(n)}=\rho^{\e_k}_s\(\big(\cL^{X^{\e_k}}\psi\big)(\cdot,Z^{\e_k}_s)\), \qquad a.s.\mP,\\
&&\lim\limits_{n\rightarrow\infty}\left(\bar{\rho}_s\big(\bar{\cL}\psi\big)\right)^{(n)}=\bar{\rho}_s\(\bar{\cL}\psi\), \qquad a.s.\mP.
\de
Moreover, by the dominated convergence theorem, it holds that
\ce
&&\lim\limits_{n\rightarrow\infty}\mE^{\mP^{\e_k}}\left(\int_0^t\left|\left(\rho^{\e_k}_s\(\big(\cL^{X^{\e_k}}\psi\big)(\cdot,Z^{\e_k}_s)\)\right)^{(n)}-\rho^{\e_k}_s\(\big(\cL^{X^{\e_k}}\psi\big)(\cdot,Z^{\e_k}_s)\)\right|^2\dif s\right)=0,\\
&&\lim\limits_{n\rightarrow\infty}\mE^{\mP^{\e_k}}\left(\int_0^t\left|\left(\bar{\rho}_s\big(\bar{\cL}\psi\big)\right)^{(n)}-\bar{\rho}_s\(\bar{\cL}\psi\)\right|^2\dif s\right)=0.
\de
So, the H\"older inequality admits us to obtain that 
\ce
&&\lim\limits_{n\rightarrow\infty}\mE^{\mP^{\e_k}}\left|\int_0^t\left(\rho^{\e_k}_s\(\big(\cL^{X^{\e_k}}\psi\big)(\cdot,Z^{\e_k}_s)\)\right)^{(n)}\dif s-\int_0^t\rho^{\e_k}_s\(\big(\cL^{X^{\e_k}}\psi\big)(\cdot,Z^{\e_k}_s)\)\dif s\right|^2=0,\\
&&\lim\limits_{n\rightarrow\infty}\mE^{\mP^{\e_k}}\left|\int_0^t\left(\bar{\rho}_s\big(\bar{\cL}\psi\big)\right)^{(n)}\dif s-\int_0^t\bar{\rho}_s\(\bar{\cL}\psi\)\dif s\right|^2=0.
\de
And then it follows from this that
\ce
&&\int_0^t\rho^{\e_k}_s\(\big(\cL^{X^{\e_k}}\psi\big)(\cdot,Z^{\e_k}_s)\)\dif s-\int_0^t\bar{\rho}_s\big(\bar{\cL}\psi\big)\dif s\\
&=&\int_0^t\rho^{\e_k}_s\(\big(\cL^{X^{\e_k}}\psi\big)(\cdot,Z^{\e_k}_s)\)\dif s-\int_0^t\left(\rho^{\e_k}_s\(\big(\cL^{X^{\e_k}}\psi\big)(\cdot,Z^{\e_k}_s)\)\right)^{(n)}\dif s\\
&&+\int_0^t\left(\rho^{\e_k}_s\(\big(\cL^{X^{\e_k}}\psi\big)(\cdot,Z^{\e_k}_s)\)\right)^{(n)}\dif s-\int_0^t\left(\bar{\rho}_s\big(\bar{\cL}\psi\big)\right)^{(n)}\dif s\\
&&+\int_0^t\left(\bar{\rho}_s\big(\bar{\cL}\psi\big)\right)^{(n)}\dif s-\int_0^t\bar{\rho}_s\big(\bar{\cL}\psi\big)\dif s\\
&=&\int_0^t\rho^{\e_k}_s\(\big(\cL^{X^{\e_k}}\psi\big)(\cdot,Z^{\e_k}_s)\)\dif s-\int_0^t\left(\rho^{\e_k}_s\(\big(\cL^{X^{\e_k}}\psi\big)(\cdot,Z^{\e_k}_s)\)\right)^{(n)}\dif s\\
&&+\sum_{j=0}^{n-1}\left(\rho^{\e_k}_{(j+1)t/n}\(\big(\cL^{X^{\e_k}}\psi\big)(\cdot,Z^{\e_k}_{(j+1)t/n})\)-\bar{\rho}_{(j+1)t/n}\big(\bar{\cL}\psi\big)\right)\((j+1)t/n-jt/n\)\\
&&+\int_0^t\left(\bar{\rho}_s\big(\bar{\cL}\psi\big)\right)^{(n)}\dif s-\int_0^t\bar{\rho}_s\big(\bar{\cL}\psi\big)\dif s\\
&\xrightarrow{w.}& 0, \qquad k\rightarrow\infty,
\de
i.e.
\be
\int_0^t\rho^{\e_k}_s\(\big(\cL^{X^{\e_k}}\psi\big)(\cdot,Z^{\e_k}_s)\)\dif s\xrightarrow{w.}\int_0^t\bar{\rho}_s\big(\bar{\cL}\psi\big)\dif s, \qquad k\rightarrow\infty.
\label{int1}
\ee

In the following, we treat the second integral in Eq.(\ref{zakai1}). By the similar deduction to above one could have that $\rho^{\e_k}_s\(\psi h(\cdot,Z^{\e_k}_s)^i\)$ converges weakly to $\bar{\rho}_s\(\psi \bar{h}^i\)$ as $k\rightarrow\infty$. Besides, define 
\ce
&&\left(\rho^{\e_k}_s\(\psi h(\cdot,Z^{\e_k}_s)^i\)\right)^{(n)}:=\sum_{j=0}^{n-1}\rho^{\e_k}_{(j+1)t/n}\(\psi h(\cdot,Z^{\e_k}_{(j+1)t/n})^i\)I_{(jt/n, (j+1)t/n]}(s),\\
&&\left(\bar{\rho}_s\(\psi \bar{h}^i\)\right)^{(n)}:=\sum_{j=0}^{n-1}\bar{\rho}_{(j+1)t/n}\(\psi \bar{h}^i\)I_{(jt/n, (j+1)t/n]}(s),
\de
and then 
\ce
&&\lim\limits_{n\rightarrow\infty}\left(\rho^{\e_k}_s\(\psi h(\cdot,Z^{\e_k}_s)^i\)\right)^{(n)}=\rho^{\e_k}_s\(\psi h(\cdot,Z^{\e_k}_s)^i\), \qquad a.s.\mP,\\
&&\lim\limits_{n\rightarrow\infty}\left(\bar{\rho}_s\(\psi \bar{h}^i\)\right)^{(n)}=\bar{\rho}_s\(\psi \bar{h}^i\), \qquad a.s.\mP.
\de
Furthermore it follows from the dominated convergence theorem that
\ce
&&\lim\limits_{k\rightarrow\infty}\mE^{\mP^{\e_k}}\left(\int_0^t\left|\left(\rho^{\e_k}_s\(\psi h(\cdot,Z^{\e_k}_s)^i\)\right)^{(n)}-\rho^{\e_k}_s\(\psi h(\cdot,Z^{\e_k}_s)^i\)\right|^2\dif s\right)=0, \\
&&\lim\limits_{k\rightarrow\infty}\mE^{\mP^{\e_k}}\left(\int_0^t\left|\left(\bar{\rho}_s\(\psi \bar{h}^i\)\right)^{(n)}-\bar{\rho}_s\(\psi \bar{h}^i\)\right|^2\dif s\right)=0.
\de
Based on the It\^o isometry, it holds that
\ce
&&\mE^{\mP^{\e_k}}\left|\int_0^t\left(\rho^{\e_k}_s\(\psi h(\cdot,Z^{\e_k}_s)^i\)\right)^{(n)}\dif \bar{B}^i_s-\int_0^t\rho^{\e_k}_s\(\psi h(\cdot,Z^{\e_k}_s)^i\)\dif \bar{B}^i_s\right|^2\\
&=&\sum_{i=1}^m\mE^{\mP^{\e_k}}\left(\int_0^t\left|\left(\rho^{\e_k}_s\(\psi h(\cdot,Z^{\e_k}_s)^i\)\right)^{(n)}-\rho^{\e_k}_s\(\psi h(\cdot,Z^{\e_k}_s)^i\)\right|^2\dif s\right),\\
&&\mE^{\mP^{\e_k}}\left|\int_0^t\left(\bar{\rho}_s\(\psi \bar{h}^i\)\right)^{(n)}\dif \bar{B}^i_s-\int_0^t\bar{\rho}_s\(\psi \bar{h}^i\)\dif \bar{B}^i_s\right|^2\\
&=&\sum_{i=1}^m\mE^{\mP^{\e_k}}\left(\int_0^t\left|\left(\bar{\rho}_s\(\psi \bar{h}^i\)\right)^{(n)}-\bar{\rho}_s\(\psi \bar{h}^i\)\right|^2\dif s\right).
\de
Thus, $\int_0^t\left(\rho^{\e_k}_s\(\psi h(\cdot,Z^{\e_k}_s)^i\)\right)^{(n)}\dif \bar{B}^i_s$ and $\int_0^t\left(\bar{\rho}_s\(\psi \bar{h}^i\)\right)^{(n)}\dif \bar{B}^i_s$ converge in mean square to\\ 
$\int_0^t\rho^{\e_k}_s\(\psi h(\cdot,Z^{\e_k}_s)^i\)\dif \bar{B}^i_s$ and $\int_0^t\bar{\rho}_s\(\psi \bar{h}^i\)\dif \bar{B}^i_s$, respectively. Let us compute 
\ce
&&\int_0^t\rho^{\e_k}_s\(\psi h(\cdot,Z^{\e_k}_s)^i\)\dif \bar{B}^i_s-\int_0^t\bar{\rho}_s\(\psi \bar{h}^i\)\dif \bar{B}^i_s\\
&=&\int_0^t\rho^{\e_k}_s\(\psi h(\cdot,Z^{\e_k}_s)^i\)\dif \bar{B}^i_s-\int_0^t\left(\rho^{\e_k}_s\(\psi h(\cdot,Z^{\e_k}_s)^i\)\right)^{(n)}\dif \bar{B}^i_s\\
&&+\int_0^t\left(\rho^{\e_k}_s\(\psi h(\cdot,Z^{\e_k}_s)^i\)\right)^{(n)}\dif \bar{B}^i_s-\int_0^t\left(\bar{\rho}_s\(\psi \bar{h}^i\)\right)^{(n)}\dif \bar{B}^i_s\\
&&+\int_0^t\left(\bar{\rho}_s\(\psi \bar{h}^i\)\right)^{(n)}\dif \bar{B}^i_s-\int_0^t\bar{\rho}_s\(\psi \bar{h}^i\)\dif \bar{B}^i_s\\
&=&\int_0^t\rho^{\e_k}_s\(\psi h(\cdot,Z^{\e_k}_s)^i\)\dif \bar{B}^i_s-\int_0^t\left(\rho^{\e_k}_s\(\psi h(\cdot,Z^{\e_k}_s)^i\)\right)^{(n)}\dif \bar{B}^i_s\\
&&+\sum_{j=0}^{n-1}\left(\rho^{\e_k}_{(j+1)t/n}\(\psi h(\cdot,Z^{\e_k}_{(j+1)t/n})^i\)-\bar{\rho}_{(j+1)t/n}\(\psi \bar{h}^i\)\right)\left(\bar{B}^i_{(j+1)t/n}-\bar{B}^i_{jt/n}\right)\\
&&+\int_0^t\left(\bar{\rho}_s\(\psi \bar{h}^i\)\right)^{(n)}\dif \bar{B}^i_s-\int_0^t\bar{\rho}_s\(\psi \bar{h}^i\)\dif \bar{B}^i_s\\
&\xrightarrow{w.}& 0, \qquad k\rightarrow\infty,
\de
that is,
\be
\int_0^t\rho^{\e_k}_s\(\psi h(\cdot,Z^{\e_k}_s)^i\)\dif \bar{B}^i_s\xrightarrow{w.}\int_0^t\bar{\rho}_s\(\psi \bar{h}^i\)\dif \bar{B}^i_s, \quad k\rightarrow\infty.
\label{int2}
\ee

For the third integral in Eq.(\ref{zakai1}), by the similar deduction to the second integral it holds that
\be
\int_0^t\int_{\mU_3}\rho^{\e_k}_s\(\psi(\lambda(s,\cdot,u)-1)\)\tilde{N}_{\lambda}(\dif s, \dif u)\xrightarrow{w.}\int_0^t\int_{\mU_3}\bar{\rho}_s\(\psi(\lambda(s,\cdot,u)-1)\)\tilde{N}_{\lambda}(\dif s, \dif u), \quad k\rightarrow\infty.
\label{int3}
\ee
Combining (\ref{int3}) with (\ref{int1}) (\ref{int2}) and taking weak limits on two hand sides of Eq.(\ref{zakai1}) as $k\rightarrow\infty$, we obtain that
\ce
\bar{\rho}_t(\psi)&=&\bar{\rho}_0(\psi)+\int_0^t\bar{\rho}_s\big(\bar{\cL}\psi\big)\dif s+\int_0^t\bar{\rho}_s\(\psi \bar{h}^i\)\dif \bar{B}^i_s\\
&&+\int_0^t\int_{\mU_3}\bar{\rho}_s\(\psi(\lambda(s,\cdot,u)-1)\)\tilde{N}_{\lambda}(\dif s, \dif u).
\de

Next, we consider $\rho^0_t(\psi)$. By the similar deduction to $\rho^{\e}_t(\psi)$,
it holds that
\ce
\rho^0_t(\psi)&=&\rho^0_0(\psi)+\int_0^t\rho^0_s\(\bar{\cL}\psi\)\dif s+\int_0^t\rho^0_s\(\psi \bar{h}^i\)\dif \bar{B}^i_s\no\\
&&+\int_0^t\int_{\mU_3}\rho^0_s\(\psi(\lambda(s,\cdot,u)-1)\)\tilde{N}_{\lambda}(\dif s, \dif u).
\de
Thus, $\bar{\rho}$ and $\rho^0$ solve the same equation
\be
\rho_t(\psi)&=&\rho_0(\psi)+\int_0^t\rho_s\(\bar{\cL}\psi\)\dif s+\int_0^t\rho_s\(\psi \bar{h}^i\)\dif \bar{B}^i_s\no\\
&&+\int_0^t\int_{\mU_3}\rho_s\(\psi(\lambda(s,\cdot,u)-1)\)\tilde{N}_{\lambda}(\dif s, \dif u).
\label{zakai3}
\ee
Besides, based on Theorem 4.2 in \cite{qd}, Eq.(\ref{zakai3}) has a unique solution. So, for $t\in[0,T]$
\ce
\bar{\rho}_t(\psi)=\rho^0_t(\psi), \qquad a.s. \mP.
\de
The proof is completed.
\end{proof}

\br
Here we can only obtain weak convergence. When trying to show convergence in probability or convergence in	mean square as that in \cite[Theorem 4.1]{ps1}, we find that it is difficult due to a jump process contained in the observation process.
\er

\section{An example}\label{exam}

In the section we give an example to explain our result. 

Consider the following slow-fast system on $\mR\times\mR$:
\be
\left\{\begin{array}{l}
\dif X^\e_t=\sin(Z^\e_t)\dif t+\sigma_1\dif V_t, \\
X^\e_0=x_0,\\
\dif Z^\e_t=\frac{1}{\e}(-Z^\e_t)\dif t+\frac{1}{\sqrt{\e}}\sigma_2\dif W_t,\\
Z^\e_0=z_0,
\end{array}
\right.
\label{exequ}
\ee
where $x_0, z_0$ are real, 
\ce
\left\{\begin{array}{l}
b_1(x,z)=\sin z, \sigma_1(x,z)=\sigma_1, f_1(x,u)=0, \\
b_2(x,z)= -z, \sigma_2(x,z)=\sigma_2, f_2(x,z,u)=0,
\end{array}
\right.
\de
and $\sigma_1, \sigma_2\neq0$ are two constants. By simple calculation, we know that $b_1(x,z), b_2(x,z), \sigma_1(x,z),$\\$ \sigma_2(x,z), f_1(x,u), f_2(x,z,u)$ satisfy {\bf Assumption 1.}. Thus, the system (\ref{exequ}) has a unique strong solution denoted by $(X^\e_t,Z^\e_t)$. 

Next take any $x\in\mR^n$ and fix it. And consider the following SDE in $\mR$:
\be\left\{\begin{array}{l}
\dif Z^x_t=-Z^x_t\dif t+\sigma_2\dif W_t,\\
Z^x_0=z_0, \qquad t\geq0.
\end{array}
\right.
\label{exfixeq}
\ee
By \cite{q1}, Eq.(\ref{exfixeq}) has a unique invariant probability measure $\bar{p}(x,dz)=\frac{1}{\sqrt{\pi\sigma^2_2}}\exp\{-\frac{z^2}{\sigma^2_2}\}\dif z$. Moreover, by some calculation we know that 
$$
\left|\int_{\mR}\varphi(z)\bar{p}(x_1; \dif z)-\int_{\mR}\varphi(z)\bar{p}(x_2; \dif z)\right|=0\leq L_{\varphi}|x_1-x_2|, \quad x_1, x_2\in\mR, \quad\varphi\in\cB(\mR),
$$
Thus, {\bf Assumption 2.} holds.

Set 
\ce
\bar{b}_1(x):=\int_{-\infty}^{+\infty}b_1(x,z)\bar{p}(x,dz)=\int_{-\infty}^{+\infty}\sin z\frac{1}{\sqrt{\pi\sigma^2_2}}\exp\{-\frac{z^2}{\sigma^2_2}\}\dif z=:\bar{b}_1, \\
(\bar{\sigma}_1\bar{\sigma}^T_1)(x):=\int_{-\infty}^{+\infty}(\sigma_1\sigma^T_1)(x,z)\bar{p}(x,\dif z)=\int_{-\infty}^{+\infty}\sigma^2_1\frac{1}{\sqrt{\pi\sigma^2_2}}\exp\{-\frac{z^2}{\sigma^2_2}\}\dif z=\sigma^2_1,
\de
and then we construct the following SDE
\ce\left\{\begin{array}{l}
\dif X^0_t=\bar{b}_1\dif t+\sigma_1\dif V_t,\\
X^0_0=x_0, \qquad\qquad 0\leq t\leq T.
\end{array}
\right.
\de
Moreover, $X^0_t=x_0+\bar{b}_1t+\sigma_1V_t$ is a unique strong solution of the above equation. Thus, by Theorem \ref{conprot}, it holds that $\{X_t^\e, t\in[0,T]\}$ converges weakly to $\{X^0_t, t\in[0,T]\}$ in $D([0,T],\mR)$.

 Take the observation process
 $$
 Y_t^{\e}=\int_0^t \arctan(X_s^{\e})\dif s+B_t+\int_0^t\int_{|u|<1}u\tilde{N}_{\lambda}(\dif s, \dif u)+\int_0^t\int_{u|\geq1}uN_{\lambda}(\dif s, \dif u),
 $$
where $h(x,z)=\arctan(x)$, $f_3(s,u)=u$, $\mU_3=\{|u|<1\}$, $\nu_3$ is a finite measure on $(\mR,\mathscr{B}(\mR))$ and $0<\lambda(t,x,u)<1$ is a constant. Then $h(x,z), f_3(s,u), \lambda(t,x,u)$ satisfy {\bf Assumption 3.-4.}.

To justify {\bf Assumption 5.}, we observe the following equation
$$
Z^\e_t=z_0+\int_0^t\frac{1}{\e}(-Z^\e_s)\dif s+\frac{1}{\sqrt{\e}}\sigma_2W_t.
$$
Set $t=\e u$, and then 
$$
Z^\e_{\e u}=z_0-\int_0^uZ^\e_{\e r}\dif r+\sigma_2 W_u.
$$
By some calculation, it holds that 

(i)  
$$
\lim_{N\rightarrow\infty}\sup\limits_{\e}\mP\{|z_0|>N\}=0,
$$

(ii) for any $U>0$ and $0\leq u_1, u_2\leq U$, there is a constant $C>0$ depending on $U, z_0, \sigma_2$ such that
$$
\mE|Z^\e_{\e u_1}-Z^\e_{\e u_2}|^4\leq C|u_1-u_2|^2, ~\mbox{for any}~ \e.
$$
So, by \cite[Theorem 4.1, Page 10]{kus}, $\{Z_{\e t}^{\e}, t\in[0,T]\}$ is tight and then {\bf Assumption 5.} is right. Moreover, by Theorem \ref{filcon}, we have that for any $\psi\in\cB(\mR^n)$, $\mE[\psi(X^{\e}_t)|\mathscr{F}_t^{Y^{\e}}]$ converges weakly to $``\mE[\psi(X^{0}_t)|\mathscr{F}_t^{Y^{\e}}]"$ as $\e\rightarrow0$.

\section{Conclusions}\label{con}

In the paper, we consider a type of nonlinear filtering problems, whose signal part is a multiscale non-Gaussian process and whose observation part is a process with jumps. Firstly, we prove that the dimension for the signal system can be reduced by a homogenized approach. Secondly, convergence of the corresponding nonlinear filtering to the homogenized filtering is shown by a weak convergence technique. Finally, we give an example to explain our result.

In the future, we will investigate the following slow-fast system on $\mR^n\times\mR^m$: for $0\leq t\leq T$,
\be\left\{\begin{array}{l}
\dif X^\e_t=b_1(X^\e_t,Z^\e_t)\dif t+\sigma_1(X^\e_t,Z^\e_t)\dif V_t+\int_{\mU_1}f_1(X^\e_{t-}, Z^\e_{t-}, u)\tilde{N}_{p_1}(\dif t, \dif u), \\
X^\e_0=x_0,\\
\dif Z^\e_t=\frac{1}{\e}b_2(X^\e_t,Z^\e_t)\dif t+\frac{1}{\sqrt{\e}}\sigma_2(X^\e_t,Z^\e_t)\dif W_t+\int_{\mU_2}f_2(X^\e_{t-},Z^\e_{t-},u)\tilde{N}^{\e}_{p_2}(\dif t, \dif u),\\
Z^\e_0=z_0.
\end{array}
\right.
\label{Eq11}
\ee
That is, the jump diffusion coefficient of the slow part contains the fast part. Although the system (\ref{Eq11}) has appeared in \cite{kus}, its dimension reduction is not essentially solved (c.f. \cite[P. 89]{kus}). Therefore, we explore the possibility of doing a similar reduction and the related nonlinear filtering.

\bigskip

\textbf{Acknowledgements:}

The author would like to thank Professor Xicheng Zhang for his valuable discussions.


\begin{thebibliography}{999}

\bibitem{Dit} P. D. Ditlevsen: Observation of $\alpha-$stable noise induced millennial climate changes from an ice record, \emph{Geophysical Research Letters}, 26(1999)1441-1444.

\bibitem{iw} N. Ikeda and S. Watanabe: {\it Stochastic differential equations
and diffusion processes,} 2nd ed., North-Holland/Kodanska,
Amsterdam/Tokyo, 1989.

\bibitem{ImkellerSri} P. Imkeller, N. S. Namachchivaya, N. Perkowski and H. C. Yeong: Dimensional reduction in
nonlinear filtering: a homogenization approach, {\it The Annals of Applied Probability,} 23(2013)2290-2326.

\bibitem{kur} T. G. Kurtz: {\it Approximation of Population Processes,} Vol. 36 of
CBMS-NSF Regional Conf. Series in Appl. Math., SIAM, Philadelphia, 1981.

\bibitem{kus} H. J. Kushner: {\it Weak Convergence Methods and Singularly Perturbed Stochastic
Control and Filtering Problems,} Systems \& Control: Foundations \& Applications 3.
Birkh\"auser, Boston, 1990.

\bibitem{kus2} H. J. Kushner: Robustness	and convergence of approximations to nonlinear filters for	jump-diffusions, {\it Computational and Applied	Math.,} 16(1996)153-183.	

\bibitem{lh} V. M. Lucic and A. J. Heunis: Convergence of nonlinear filters for randomly
perturbed dynamical systems, {\it Appl. Math. Optim.,} 48(2003)93-128.

\bibitem{mbfp} T. Meyer-Brandis and F. Proske: Explicit solution of a non-linear filtering problems for
L\'evy Processes with application to finance, {\it Applied Mathematics and Optimization,} 50(2004)119-134.

\bibitem{ps1} A. Papanicolaou and K. Spiliopoulos: Dimension reduction in statistical estimation of partially observed multiscale processes, {\it SIAM J. on Uncertainty Quantification,} 5(2017)1220-1247.

\bibitem{ps2} A. Papanicolaou and K. Spiliopoulos: Filtering the	maximum	likelihood	for	
multiscale	problems, {\it SIAM J. on Multiscale Model. Simul.,} 12(2014)1193-1229.

\bibitem{ap} A. Pazy: {\it Semigroups of Linear Operators and Applications to
Partial Differential Equations,} Springer-Verlag Berlin Heidelberg New York, 1983.

\bibitem{pny}  J. H. Park, N. S. Namachchivaya and H. C. Yeong: Particle filters in a multiscale environment: Homogenized hybrid particle filter,{\it J. Appl. Mech.,} 78(2011)1-10. 

\bibitem{psn}  J. H. Park, R. B. Sowers and N. S. Namachchivaya: Dimensional reduction
in nonlinear filtering, {\it Nonlinearity,} 23(2010)305-324.

\bibitem{q0} H. J. Qiao: Effective filtering for multiscale stochastic dynamical systems
driven by L\'evy processes, https://arxiv.org/abs/1810.10370.

\bibitem{q1} H. J. Qiao: Exponential ergodicity for SDEs with jumps and non-Lipschitz coefficients,
{\it Journal of Theoretical Probability}, 27(2014)137-152.

\bibitem{q2} H. J. Qiao: Euler-Maruyama approximation for SDEs with jumps and non-Lipschitz coefficients,
{\it Osaka Journal of Mathematics}, 51(2014)47-66.

\bibitem{q3} H. J. Qiao: Uniqueness for measure-valued equations of nonlinear filtering     for stochastic dynamical systems with L\'evy noises, {\it Advances in Applied Probability,} 50(2018)396-413.

\bibitem{qd} H. J. Qiao and J. Q. Duan: Nonlinear Filtering of Stochastic Dynamical Systems with L\'evy Noises,
{\it Advances in Applied Probability,} 47(2015)902-918.

\bibitem{st} D. W. Stroock: Diffusion processes associated with L\'evy generators. {\it Z. Wahrscheinlichkeitstheorie und Verw. Gebiete,} 32 (1975)209-244.

\bibitem{sv} D. W. Stroock and S. R. S. Varadhan: Diffusion processes with boundary conditions, {\it Comm. Pure Appl. Math.,}
24(1971)147-225.

\bibitem{yz} G. Yin and Q. Zhang: {\it Continuous Time Markov Chains and Applications: A Two Time-Scale Approach,} Springer-Verlag, New York, 2013.	

\bibitem{zqd} Y. J. Zhang, H. J. Qiao and J. Q. Duan : Effective filtering analysis for non-Gaussian dynamic systems. Appear in {\it Applied Mathematics and Optimization}.
\end{thebibliography}
\end{document}